\documentclass[11pt]{article}
\usepackage{amssymb}

\usepackage{amsmath}
\usepackage{theorem}
\usepackage{here}
\usepackage{color}

\newtheorem{thm}{Theorem}[section]
\newtheorem{lem}{Lemma}[section]
\newtheorem{cor}{Corollary}[section]
\newtheorem{prp}{Proposition}[section]
\theorembodyfont{\rmfamily}

\newtheorem{remark}{Remark}[section]
\makeatletter

\@addtoreset{equation}{section}
\makeatother

\usepackage{comment} %
\usepackage{bm}
\usepackage[dvipdfmx]{graphicx} %

\def\al{{\alpha}}
\def\be{{\beta}}

\def\de{{\delta}}
\def\ep{{\varepsilon}}
\def\la{{\lambda}}

\def\om{{\omega}}
\def\th{{\theta}}

\def\ka{{\kappa}}

\def\bla{{\text{\boldmath $\lambda$}}}

\def\bth{{\text{\boldmath $\theta$}}}

\def\bxi{{\text{\boldmath $\xi$}}}

\def\De{{\Delta}}

\def\Ga{{\Gamma}}

\def\La{{\Lambda}}
\def\a{{\text{\boldmath $a$}}}

\def\d{{\text{\boldmath $d$}}}
\def\e{{\text{\boldmath $e$}}}

\def\j{{\text{\boldmath $j$}}}

\def\p{{\text{\boldmath $p$}}}

\def\v{{\text{\boldmath $v$}}}

\def\x{{\text{\boldmath $x$}}}

\def\z{{\text{\boldmath $z$}}}

\def\D{{\text{\boldmath $D$}}}

\def\X{{\text{\boldmath $X$}}}

\def\ah{{\hat a}}

\def\ph{{\hat p}}

\def\diag{{\rm diag\,}}

\def\[{{\text{\boldmath $[$}}}
\def\]{{\text{\boldmath $]$}}}

\def\/{{\Bigr/\!\!}}

\def\1r{{\rm (1)}}
\def\2r{{\rm (2)}}
\def\3r{{\rm (3)}}
\def\4r{{\rm (4)}}
\def\5r{{\rm (5)}}

\def\non{{\nonumber}}
\begin{document}
\title{Bayesian Shrinkage Estimation of %
Negative Multinomial Parameter Vectors}
\author{
Yasuyuki Hamura\footnote{Graduate School of Economics, University of Tokyo, 
7-3-1 Hongo, Bunkyo-ku, Tokyo 113-0033, JAPAN.\newline{
E-Mail: yasu.stat@gmail.com}} \
and
Tatsuya Kubokawa\footnote{Faculty of Economics, University of Tokyo, 
7-3-1 Hongo, Bunkyo-ku, Tokyo 113-0033, JAPAN. \newline{
E-Mail: tatsuya@e.u-tokyo.ac.jp}}
}
\maketitle
\begin{abstract}
The negative multinomial distribution is a multivariate generalization of the negative binomial distribution. 
In this paper, we consider the problem of estimating an unknown matrix of probabilities on the basis of observations of negative multinomial variables under the standardized squared error loss. 
First, a general sufficient condition for a shrinkage estimator to dominate the UMVU estimator is derived and an empirical Bayes estimator satisfying the condition is constructed. 
Next, a hierarchical shrinkage prior is introduced, an associated Bayes estimator is shown to dominate the UMVU estimator under some conditions, and some remarks about posterior computation are presented. 
Finally, shrinkage estimators and the UMVU estimator are compared by simulation. 
\par\vspace{4mm}
\end{abstract}

\section{Introduction}
\label{sec:introduction}
Stein's phenomenon for the estimation of parameters of discrete distributions has been extensively studied since Clevenson and Zidek (1975) showed that the usual estimator of the mean vector of independent Poisson distributions is dominated by a Bayesian shrinkage estimator under the standardized squared error loss. 
For example, Ghosh and Parsian (1981), Tsui (1979b), Tsui and Press (1982), and Ghosh and Yang (1988) considered different estimators of Poisson parameters under different loss functions. 
Estimation for discrete exponential families including the Poisson and the negative binomial distributions was treated by Tsui (1979a), Hwang (1982), and Ghosh, Hwang, and Tsui (1983). 
Tsui (1984), Tsui (1986a), and Tsui (1986b) explored the robustness of Clevenson-Zidek-type estimators in estimating means when the observations are not Poisson-distributed. 
In particular, Tsui (1986b) considered the case of dependent observations following the negative multinomial distribution, which is a multivariate generalization of the negative binomial distribution and arises as the joint distribution of the frequencies of multiple events in inverse sampling. 
The negative multinomial distribution is also included in the general classes of discrete distributions of Chou (1991) and Dey and Chung (1992). 

However, little attention has been paid to the construction of Bayesian shrinkage estimators when the underlying distributions are not Poisson. 
This could be partly because tractable hierarchical models may not be so widely known in such cases; some difficulties with the beta-binomial hierarchy are discussed in Example 4.5.3 of Lehmann and Casella (1998). 
In this paper, we consider the Bayesian estimation of multiple negative multinomial parameter vectors. 

The $m$-dimensional negative multinomial distribution with parameters $r > 0$ and $\mathring{\p } = ( \mathring{p} _1 , \dots , \mathring{p} _m )' \in D_m = \big\{ ( \tilde{p} _1 , \dots , \tilde{p} _m )' | \tilde{p} _1 , \dots , \tilde{p} _m > 0, \, \sum_{i = 1}^{m} \tilde{p} _i < 1\big\} $, denoted by ${\rm{NM}}_m (r, \mathring{\p } )$, has probability mass function 
\begin{align}
{\rm{NM}}_m ( \x | r, \mathring{\p } ) &= {\Ga \big( r + \sum_{i = 1}^{m} x_i \big) \over \Ga (r) \prod _{i = 1}^{m} x_i !} {\mathring{p} _0}^r \prod_{i = 1}^{m} {\mathring{p} _i}^{x_i} \label{eq:NM} 
\end{align}
for $\x = ( x_1 , \dots , x_m )' \in {\mathbb{N} _0}^m = \{ 0, 1, 2, \dotsc \} ^m$, where %
$\mathring{p} _0 = 1 - \mathring{p} _{\cdot } = 1 - \sum_{i = 1}^{m} \mathring{p} _i$ and where $r$ corresponds to the number of successes in inverse sampling. 
Even if $r$ is not an integer, the probability function (\ref{eq:NM}) is well defined and has the Poisson-gamma mixture representation 
\begin{align}
{\rm{NM}}_m ( \x | r, \mathring{\p } ) &= \int_{0}^{\infty } {v^{r - 1} \over \Ga (r)} e^{- v} \Big[ \prod_{i = 1}^{m} {\{ ( \mathring{p} _i / \mathring{p} _0 ) v \} ^{x_i} \over {x_i}!} e^{- ( \mathring{p} _i / \mathring{p} _0 ) v} \Big] dv \text{.} \label{eq:NM_mixture} 
\end{align}
The mean and variance of the negative multinomial distribution ${\rm{NM}}_m (r, \mathring{\p } )$ are $r \mathring{\p } / \mathring{p} _0$ and $r \diag ( \mathring{\p } ) / \mathring{p} _0 + r \mathring{\p } {\mathring{\p }}' / {\mathring{p} _0}^2$. 
The marginals are negative binomial. 
If $\mathring{\X } ^{(1)} \sim {\rm{NM}}_m (r^{(1)} , \mathring{\p } )$ and $\mathring{\X } ^{(2)} \sim {\rm{NM}}_m ( r^{(2)} , \mathring{\p } )$ for %
$r^{(1)} , r^{(2)} > 0$, then $\mathring{\X } ^{(1)} + \mathring{\X } ^{(2)} \sim {\rm{NM}}_m (r^{(1)} + r^{(2)} , \mathring{\p } )$; therefore, $r$ can also be interpreted as a sample size. 
For further properties and applications of the negative multinomial distribution, see, for example, Sibuya, Yoshimura, and Shimizu (1964) and Tsui (1986b) and the references therein. 

Suppose that $\X _1 = ( X_{1, 1} , \dots , X_{m, 1} )' , \dots , \X _N = ( X_{1, N} , \dots , X_{m, N} )'$ are independently distributed as ${\rm{NM}}_m (r, \p _1 ), \dots , {\rm{NM}}_m (r, \p _N )$, respectively, for $m, N \in \mathbb{N} = \{ 1, 2, \dotsc \} $, where all the elements of $\p = ( \p _1 , \dots , \p _N ) = ( ( p_{1, 1} , \dots , p_{m, 1} )' , \dots , ( p_{1, N} , \dots , p_{m, N} )' ) \in {D_m}^N$ are assumed to be unknown. 
For $n = 1, \dots , N$, we consider the problem of estimating %
the matrix $( \p _1 , \dots , \p _n )$ on the basis of $\X = ( \X _1 , \dots , \X _N )$ under the standardized squared error loss 
\begin{align}
L_n ( \d , \p ) = \sum_{\nu = 1}^{n} \sum_{i = 1}^{m} {1 \over p_{i, \nu }} ( d_{i, \nu } - p_{i, \nu } )^2 \text{,} \label{eq:loss_SS} 
\end{align}
where $\d = ( d_{i, \nu } )_{1 \le i \le m, \, 1 \le \nu \le N} \in \mathbb{R} ^{m \times N}$. 
Here, $n = N$ corresponds to the simultaneous estimation of all the parameters while $n = 1$ corresponds to the estimation of $\p _1$ relating to the first observation $\X _1$ by using all the information $\X $. 

As prior distribution for $\p $, we first use the conjugate Dirichlet distribution with density 
\begin{align}
\prod_{\nu = 1}^{N} {\rm{Dir}} _m ( \p _{\nu } | a_0 , \a ) = \prod_{\nu = 1}^{N} \Big\{ {\Ga ( a_0 + a _{\cdot } ) \over \Ga ( a_0 ) \prod_{i = 1}^{m} \Ga ( a_i )} {p_{0, \nu }}^{a_0 - 1} \prod_{i = 1}^{m} {p_{i, \nu }}^{a_i - 1} \Big\} \text{,} \label{eq:Dir} 
\end{align}
where $a_0 \in \mathbb{R}$, $\a = ( a_1 , \dots , a_m )' \in (0, \infty )^m$, $a_{\cdot } = \sum_{i = 1}^{m} a_i$, and $p_{0, \nu } = 1 - p_{\cdot , \nu } = 1 - \sum_{i = 1}^{m} p_{i, \nu }$ for $\nu = 1, \dots , N$. %
As will be shown later, the UMVU estimator of $\p $ is $\hat{\p } ^{\rm{U}} = ( X_{i, \nu } / (r + X_{\cdot , \nu } - 1))_{1 \le i \le m, \, 1 \le \nu \le N}$, where $X_{\cdot , \nu } = \sum_{i = 1}^{m} X_{i, \nu }$ for $\nu = 1, \dots , N$, and corresponds to the Bayes estimator with respect to the prior (\ref{eq:Dir}) with $a_0 = - m$ and $\a = \j ^{(m)}$ and the loss (\ref{eq:loss_SS}) with $n = N$, where $\j ^{(m)} = (1, \dots , 1)' \in \mathbb{R} ^m$. 
Also, it will be seen that the Jeffreys prior is (\ref{eq:Dir}) with $a_0 = - (m - 1) / 2$ and $\a = \j ^{(m)} / 2$. 

In Section \ref{sec:eb}, we first consider the general class of estimators 
\begin{align}
\hat {\p } ^{( \de )} = \Big( {X_{i, \nu } \over r + X_{\cdot , \nu } - 1 + \de ( X_{\cdot , \cdot } )} \Big) _{1 \le i \le m, \, 1 \le \nu \le N} \text{,} \quad \label{eq:shrinkage_estimator_general} 
\end{align}
where $\de ( X_{\cdot , \cdot } )$ is a strictly positive function of $X_{\cdot , \cdot } = \sum_{\nu = 1}^{N} X_{\cdot , \nu } = \sum_{\nu = 1}^{N} \sum_{i = 1}^{m} X_{i, \nu }$, and derive a sufficient condition for the shrinkage estimator $\hat {\p } ^{( \de )}$ to dominate the unbiased estimator $\hat{\p } ^{\rm{U}}$. 
Next we construct an empirical Bayes estimator based on the prior (\ref{eq:Dir}) with $\a = \j ^{(m)}$ and show that it dominates the unbiased estimator when $m$ is sufficiently large by using the derived condition. 

In Section \ref{sec:hb}, we obtain a shrinkage estimator of the form $( X_{i, \nu } / \{ r + X_{\cdot , \nu } - 1 + \de ( \X _{\cdot } ) \} )_{1 \le i \le m, \, 1 \le \nu \le N}$, where $\de ( \X _{\cdot } ) > 0$ is some symmetric function of $\X _{\cdot } = ( X_{\cdot , 1} , \dots , X_{\cdot , N} )'$, by introducing a hierarchical prior for $\p $. 
In a simple case, this prior becomes 
\begin{align}
\p \sim \Big( \prod_{\nu = 1}^{N} p_{0, \nu } \Big) ^{- m - 1} / \Big( \sum_{\nu = 1}^{N} \log {1 \over p_{0, \nu }} \Big) ^{\al } \text{,} \non 
\end{align}
where $\al > 0$. 
The above expression shows that the prior puts more probability around $p_{0, 1} = \dots = p_{0, N} = 1$ than the Dirichlet prior $\p \sim \prod_{\nu = 1}^{N} {p_{0, \nu }}^{- m - 1}$. 
Our hierarchical Bayes estimator is shown to dominate the UMVU estimator under some conditions. 
Also, for sufficiently large $m$, we obtain an estimator based on our hierarchical prior which dominates a Bayes estimator against the Jeffreys prior under the loss 
\begin{align}
\tilde{L} _n ( \tilde{\d } , \p ) = \sum_{\nu = 1}^{n} \sum_{i = 1}^{m} \Big( \tilde{d} _{i, \nu } - p_{i, \nu } - p_{i, \nu } \log {\tilde{d} _{i, \nu } \over p_{i, \nu }} \Big) \text{,} \label{eq:loss_KL} 
\end{align}
where $\tilde{\d } = ( \tilde{d} _{i, \nu } )_{1 \le i \le m, \, 1 \le \nu \le N} \in (0, \infty )^{m \times N}$. 
In addition, it turns out that posterior computation is quite simple under our hierarchical prior. 

Recently, Stoltenberg and Hjort (2019) also considered Bayesian multivariate models for count variables based on the Poisson likelihood. 
Hamura and Kubokawa (2019) and Hamura and Kubokawa (2020) considered estimation of Poisson parameters when sample sizes are unbalanced by using and generalizing the shrinkage prior of Komaki (2015). 
Interestingly, it is the method for evaluating integrals in Bayesian predictive probabilities of Poisson variables in the presence of unbalanced sample sizes, developed by Komaki (2015) and utilized by Hamura and Kubokawa (2019) and Hamura and Kubokawa (2020), that plays a crucial role in obtaining the results in Section \ref{sec:hb} for our hierarchical Bayes estimators of negative multinomial parameters in the balanced setting. 

The remainder of the paper is organized as follows. 
In Sections \ref{sec:eb} and \ref{sec:hb}, we consider empirical Bayes and hierarchical Bayes estimators, respectively. 
In Section \ref{sec:sim}, through simulation, we compare our proposed estimators with the UMVU estimator as well as an alternative estimator which estimates $\p _1 , \dots , \p _N$ independently based on $\X _1 , \dots , \X _N$, respectively. 
Proofs are in the Appendix.

\section{Empirical Bayes Estimation}
\label{sec:eb}
We first derive a sufficient condition for the shrinkage estimator $\hat{\p } ^{( \de )}$ given in (\ref{eq:shrinkage_estimator_general}) to dominate the UMVU estimator. 
Let 
\begin{align}
\ph _{i, \nu }^{\rm{U}} = \begin{cases} \displaystyle {X_{i, \nu } \over r + X_{\cdot , \nu } - 1} & \text{if $X_{i, \nu } \ge 1$} \\ \displaystyle 0 & \text{if $X_{i, \nu } = 0$} \end{cases} \label{eq:UMVU} 
\end{align}
for $i = 1, \dots , m$ and $\nu = 1, \dots , N$. 
Then $\hat{\p } ^{\rm{U}} = ( \ph _{i, \nu }^{\rm{U}} )_{1 \le i \le m, \, 1 \le \nu \le N}$ is the UMVU estimator of $\p $ since it is unbiased by Lemma \ref{lem:hudson} in the Appendix and since $\x $ is a complete sufficient statistic. 
Let $\de \colon \mathbb{N} _0 \to (0, \infty )$, 
\begin{align}
\ph _{i, \nu }^{( \de )} &= \begin{cases} \displaystyle {X_{i, \nu } \over r + X_{\cdot , \nu } - 1 + \de ( X_{\cdot , \cdot } )} & \text{if $X_{i, \nu } \ge 1$} \\ \displaystyle 0 & \text{if $X_{i, \nu } = 0$} \end{cases} \non 
\end{align}
for $i = 1, \dots , m$ and $\nu = 1, \dots , N$, and $\hat{\p } ^{( \de )} = ( \ph _{i, \nu }^{( \de )} )_{1 \le i \le m, \, 1 \le \nu \le N}$. 

\begin{thm}
\label{thm:shrinkage_estimator_general} 
Let $n = 1, \dots , N$ and assume $r \ge 5 / 2$. 
Suppose that the function $\de $ satisfies the following conditions for all $z \in \mathbb{N}$: 
\begin{enumerate}
\item
$z \de (z) \le (z + 1) \de (z + 1)$. 
\item
If $z \ge 2$, then 
\begin{itemize}
\item
$\de (z) \le 2 (m - 3)$ implies $(m - 6) \de (z) + 2 (m - 3) r \ge 0$ and 
\item
$\de (z) > 2 (m - 3)$ implies $n \{ (m - 6) \de (z) + 2 (m - 3) r \} \ge (z - 1) \{ \de (z) - 2 (m - 3) \} $. 
\end{itemize}
\end{enumerate}
Then the shrinkage estimator $\hat{\p } ^{( \de )}$ %
dominates the UMVU estimator $\hat{\p } ^{\rm{U}}$ under the loss $L_n ( \d , \p )$ given by (\ref{eq:loss_SS}). 
\end{thm}

\noindent
For example, if $\de ( X_{\cdot , \cdot } ) = c_0$ for some constant $0 < c_0 \le 2 (m - 3)$, conditions (i) and (ii) are satisfied provided that $m \ge 6 (r + c_0 ) / (2 r + c_0 )$. 
Also, condition (i) is satisfied if $\de ( X_{\cdot , \cdot } ) = c_1 + c_2 / X_{\cdot , \cdot }$ for some constants $c_1 , c_2 > 0$ when $X_{\cdot , \cdot } \ge 1$. 

Next, we construct an empirical Bayes estimator. 
Lemma \ref{lem:B} below states that the shrinkage estimator $\hat{\p } ^{( \de )}$ coincides with a Bayes solution in a simple case. 
Let $\de ^{( a_0 )} ( X_{\cdot , \cdot } ) = a_0 + m$. 

\begin{lem}
\label{lem:B}
Suppose $a_0 > \max \{ - m, - r \} $. 
Then the shrinkage estimator $\hat{\p } ^{( \de ^{( a_0 )} )}$ 
is a Bayes solution with respect to the prior (\ref{eq:Dir}) with $\a = \j ^{(m)}$ under the loss (\ref{eq:loss_SS}) for every $n = 1, \dots , N$. 
\end{lem}

\noindent
The conditions $a_0 > - m$ and $a_0 > - r$ ensure, respectively, %
that $\hat{\p } ^{( \de ^{( a_0 )} )}$ shrinks toward the origin and that the posterior distribution is proper. 
Additionally, it follows from Theorem \ref{thm:shrinkage_estimator_general} and Lemma \ref{lem:B} that if $r \ge 5 / 2$ and $0 < a_0 \le m - 6$, the estimator $\hat{\p } ^{( \de ^{( a_0 )} )} = ( X_{i, \nu } / (r + a_0 + X_{\cdot , \nu } + m - 1))_{1 \le i \le m, \, 1 \le \nu \le N}$ is proper Bayes and dominates the UMVU estimator. 

An empirical Bayes estimator is obtained by first assuming $a_0 > 1$ and then substituting for $a_0$ in $\hat{\p } ^{( \de ^{( a_0 )} )}$ an estimator %
based on the marginal likelihood of $\p $ under the prior corresponding to $\hat{\p } ^{( \de ^{( a_0 )} )}$. 
More specifically, when $a_0 > 1$, the prior expectation of the mean $E[ X_{\cdot , \cdot } ] = \sum_{\nu = 1}^{N} \sum_{i = 1}^{m} r p_{i, \nu } / p_{0, \nu }$ with respect to the Dirichlet prior (\ref{eq:Dir}) with $\a = \j ^{(m)}$ is given by 
\begin{align}
\int_{{D_m}^N} E[ X_{\cdot , \cdot } ] \Big\{ \prod_{{\nu }' = 1}^{N} {\rm{Dir}}_m ( \p _{{\nu }'} | a_0 , \j ^{(m)} ) \Big\} d\p &= \sum_{\nu = 1}^{N} \sum_{i = 1}^{m} r \int_{{D_m}^N} {p_{i, \nu } \over p_{0, \nu }} \Big\{ \prod_{{\nu }' = 1}^{N} {\rm{Dir}}_m ( \p _{{\nu }'} | a_0 , \j ^{(m)} ) \Big\} d\p \non \\
&= {N m r \over a_0 - 1} \text{.} \non 
\end{align}
Thus, an estimator of $a_0$ is obtained as 
\begin{align}
\ah _0 &= 1 + N m r / X_{\cdot , \cdot } \non 
\end{align}
and our empirical Bayes estimator is 
\begin{align}
\hat{\p } ^{\rm{EB}} &= ( \ph _{i, \nu }^{\rm{EB}} )_{1 \le i \le m, \, 1 \le \nu \le N} = \hat{\p } ^{( \de ^{( a_0 )} )} |_{a_0 = \ah _0} \non \\
&= \Big( {X_{i, \nu } \over r + X_{\cdot , \nu } - 1 + \de ^{\rm{EB}} ( X_{\cdot , \cdot } )} \Big) _{1 \le i \le m, \, 1 \le \nu \le N} \text{,} \label{eq:EB} 
\end{align}
where 
\begin{align}
\de ^{\rm{EB}} ( X_{\cdot , \cdot } ) &= 1 + m + N m r / X_{\cdot , \cdot } \non 
\end{align}
when $X_{\cdot , \cdot } \ge 1$ and $\de ^{\rm{EB}} (0) \in (1 + m + N m r, \infty )$. 

The following corollary gives a sufficient condition for $\hat{\p } ^{\rm{EB}}$ to dominate the UMVU estimator. 

\begin{cor}
\label{cor:EB}
Suppose that $m \ge 7$ and that $r \ge 5 / 2$. 
Then $\hat{\p } ^{\rm{EB}}$ is an empirical Bayes estimator dominating the UMVU estimator $\hat{\p } ^{\rm{U}}$ under the loss $L_n ( \d , \p )$ given by (\ref{eq:loss_SS}) for every $n = 1, \dots , N$. 
\end{cor}

\noindent
It is worth noting that the condition given in the above corollary is independent of $n$, which shows some robustness of the empirical Bayes estimator $\hat{\p } ^{\rm{EB}}$. 
Additionally, we do not have to set $N > m, r$, nor do we need to assume $r > m$. 

The UMVU estimator corresponds to $a_0 = - m$ since $\lim_{a_0 \to - m} \hat{\p } ^{( \de ^{( a_0 )} )} = \hat{\p } ^{\rm{U}}$ (when $r > m$). 
However, the empirical Bayes estimator $\hat{\p } ^{\rm{EB}}$ was derived under the assumption that $a_0 > 1$. 
Indeed, we have $\ah _0 > 1$ since all the elements of the observations $\X = ( \X _1 , \dots , \X _N )$ are nonnegative. 
Thus, there is a discrepancy in the support of $a_0$ between the usual Bayes estimator and the empirical Bayes estimator. 
On the other hand, in the case of hierarchical Bayes estimation, a mixture of the priors $\p \sim \prod_{\nu = 1}^{N} {\rm{NM}}_m ( \p _{\nu } | s, \j ^{(m)} )$, $s > - m$, will be considered in the next section.

\section{Hierarchical Bayes Estimation}
\label{sec:hb}
In this section, we first introduce a shrinkage prior for $\p $ and investigate its properties (Section \ref{subsec:hb_prior}). 
Next, using the prior, we construct a hierarchical Bayes estimator that dominates the UMVU estimator under some conditions (Section \ref{subsec:hb_dominance}). 
Finally, some remarks about posterior computation are presented (Section \ref{subsec:hb_computation}).

\subsection{A hierarchical shrinkage prior}
\label{subsec:hb_prior} 
For $\p = %
( ( p_{1, 1} , \dots , p_{m, 1} )' , \dots , ( p_{1, N} , \dots , p_{m, N} )' ) \in {D_m}^N$ and $p_{0, \nu } = 1 - p_{\cdot , \nu } = 1 - \sum_{i = 1}^{m} p_{i, \nu }$, $\nu = 1, \dots , N$, let 
\begin{align}
\pi _{\al , \be , g, a_0 , \a } ( \p ) &= \int_{0}^{\infty } t^{\al - 1} e^{- \be t} g(t) \Big\{ \prod_{\nu = 1}^{N} \Big( {p_{0, \nu }}^{t + a_{0} - 1} \prod_{i = 1}^{m} {p_{i, \nu }}^{a_{i} - 1} \Big) \Big\} dt \text{,} \label{eq:shrinkage_prior} 
\end{align}
where $\al > 0$, $\be \ge 0$, %
$g \colon (0, \infty ) \to (0, \infty )$ is a bounded and smooth function, $a_0 \in \mathbb{R}$, and $\a = ( a_1 , \dots , a_m )' \in (0, \infty )^m$. 
When $g = g_1$, where $g_1 \colon (0, \infty ) \to (0, \infty )$ is the function defined by $g_1 (t) = 1$, $t \in (0, \infty )$, 
the prior (\ref{eq:shrinkage_prior}) becomes 
\begin{align}
\pi _{\al , \be , g_1 , a_0 , \a } ( \p ) = \Ga ( \al ) \Big\{ \prod_{\nu = 1}^{N} \Big( {p_{0, \nu }}^{a_{0} - 1} \prod_{i = 1}^{m} {p_{i, \nu }}^{a_{i} - 1} \Big) \Big\} / \Big( \be + \sum_{\nu = 1}^{N} \log {1 \over p_{0, \nu }} \Big) ^{\al } \text{.} \label{eq:shrinkage_prior_simple} 
\end{align}
It can be seen that 
\begin{align}
\lim_{\al \to 0} {\pi _{\al , \be , g_1 , a_0 , \a } ( \p ) \over \Ga ( \al )} = \prod_{\nu = 1}^{N} \Big( {p_{0, \nu }}^{a_{0} - 1} \prod_{i = 1}^{m} {p_{i, \nu }}^{a_{i} - 1} \Big) \propto \prod_{\nu = 1}^{N} {\rm{Dir}}_m ( \p _{\nu } | a_0 , \a ) \text{.} \non 
\end{align}
and that the denominator of (\ref{eq:shrinkage_prior_simple}) tends to infinity as $\min \{ p_{0, 1} , \dots , p_{0, N} \} \to 0$. 
Thus, $\pi _{\al , \be , g, a_0 , \a } ( \p )$ is a shrinkage prior based on the Dirichlet distribution. 
Furthermore, if $m = 1$, $N \ge 2$, and $( \La , \bth ) \sim e^{- ( a_0 - 1) \La }$, where $\La = \sum_{{\nu }' = 1}^{N} \log (1 / p_{0, {\nu }'} )$ and $\th _{\nu } = \{ \log (1 / p_{0, \nu } ) \} / \sum_{{\nu }' = 1}^{N} \log (1 / p_{0, {\nu }'} )$, $\nu = 1, \dots , N - 1$, then $\p \sim \big( \prod_{\nu = 1}^{N} {p_{0, \nu }}^{a_0 - 1} \big) / \big\{ \sum_{\nu = 1}^{N} \log (1 / p_{0, \nu } ) \big\} ^{N - 1} \propto \pi _{N - 1, 0, g_1 , a_0 , 1} ( \p )$. 

Let $a_{\cdot } = \sum_{i = 1}^{m} a_i$. 
Necessary and sufficient conditions for propriety of the prior and posterior distributions are as follows. 

\begin{lem}
\label{lem:propriety} 
\begin{enumerate}
\item
The prior (\ref{eq:shrinkage_prior}) is proper if and only if either 
\begin{itemize}
\item
$a_0 > 0$ and $\int_{1}^{\infty } t^{\al - N a_{\cdot } - 1} e^{- \be t} g(t) dt < \infty $ or 
\item
$a_0 = 0$, $\int_{0}^{1} t^{\al - N - 1} e^{- \be t} g(t) dt < \infty $, and $\int_{1}^{\infty } t^{\al - N a_{\cdot } - 1} e^{- \be t} g(t) dt < \infty $. 
\end{itemize}
\item
Under the prior (\ref{eq:shrinkage_prior}), the posterior distribution of $\p $ given the observations $\X = ( \x _1 , \dots , \x _N )$ is proper for all $\x _1 , \dots , \x _N \in {\mathbb{N} _0}^m$ if and only if either 
\begin{itemize}
\item
$r + a_0 > 0$ and $\int_{1}^{\infty } t^{\al - N a_{\cdot } - 1} e^{- \be t} g(t) dt < \infty $ or 
\item
$r + a_0 = 0$, $\int_{0}^{1} t^{\al - N - 1} e^{- \be t} g(t) dt < \infty $, and $\int_{1}^{\infty } t^{\al - N a_{\cdot } - 1} e^{- \be t} g(t) dt < \infty $. 
\end{itemize}
\end{enumerate}
\end{lem}

\noindent
When the condition of part (ii) of Lemma \ref{lem:propriety} is satisfied, we will simply say that the posterior is proper. 
For example, when $g = g_1$ and either $\al < N a_{\cdot }$ or $\be > 0$, the prior (\ref{eq:shrinkage_prior}) is proper if $a_0 > 0$, while the posterior is proper if $r + a_0 > 0$. 
It is also worth noting that even when $a_0 < 0$ and the prior is improper, the condition for posterior propriety may still be satisfied. 

The prior (\ref{eq:shrinkage_prior}) is related to shrinkage priors in the Poisson case. 
Specifically, if $m = 1$ and $\bla = ( \la _1 , \dots , \la _N ) = ( \log (1 / p_{0, 1} ), \dots , \log (1 / p_{0, N} )) \approx {\bm{0} ^{(N)}}'$, where $\bm{0} ^{(N)} = (0, \dots , 0)' \in \mathbb{R} ^N$, then $\bla $ is approximately distributed as 
\begin{align}
\bla &\sim \int_{0}^{\infty } t^{\al - 1} e^{- \be t} g(t) \Big\{ \prod_{\nu = 1}^{N} (e^{- \la _{\nu }} )^{t + a_{0}} (1 - e^{- \la _{\nu }} )^{a_{1} - 1} \Big\} dt \non \\
&\approx \Big( \prod_{\nu = 1}^{N} {\la _{\nu }}^{a_{1} - 1} \Big) \int_{0}^{\infty } t^{\al - 1} e^{- t ( \be + \la _{\cdot } )} g(t) dt \text{,} \label{eq:shrinkage_prior_limit} 
\end{align}
where $\la _{\cdot } = \sum_{\nu = 1}^{N} \la _{\nu }$. 
The density (\ref{eq:shrinkage_prior_limit}) corresponds to the prior considered by Ghosh and Parsian (1981) when $a_1 = 1$, to that considered by Komaki (2004) when $\be = 0$ and $g = g_1$, and to that considered by Komaki (2006) when $\al = m a_1 - 1$, $\be = 0$, and $g(t) = \{ t / (1 + \ka t) \} ^{c + 1}$ for all $t \in (0, \infty )$ for some $c > - m a_1$ and $\ka > 0$. 
However, in order to prove the results in the next subsection, we need to extend the technique of the proof of Theorem 1 of Komaki (2015), who considered an unbalanced problem.

\subsection{Dominance results}
\label{subsec:hb_dominance} 
In order to derive an explicit form of a Bayes solution with respect to the prior (\ref{eq:shrinkage_prior}), we define 
\begin{align}
K( \al , \be , g, \xi _0 , \bxi ) &= \int_{0}^{\infty } t^{\al - 1} e^{- \be t} g(t) \Big\{ \prod_{\nu = 1}^{N} {\Ga (t + \xi _0 ) \over \Ga (t + \xi _0 + \xi _{\nu } )} \Big\} dt \label{eq:K} 
\end{align}
for $\xi _0 \ge 0$ and $\bxi = ( \xi _1 , \dots , \xi _N )' \in [0, \infty )^N$ and we let $\j ^{(N)} = (1, \dots , 1)' \in \mathbb{R} ^N$. 
For now, we consider the case of $a_0 = - m$ and $\a = \j ^{(m)}$%
and assume that either 
\begin{align}
r > m \quad \text{and} \quad \int_{1}^{\infty } t^{\al - N m - 1} e^{- \be t} g(t) dt < \infty \label{eq:assumption1} 
\end{align}
or 
\begin{align}
r = m \text{,} \quad \int_{0}^{1} t^{\al - N - 1} e^{- \be t} g(t) dt < \infty \text{,} \quad \text{and} \quad \int_{1}^{\infty } t^{\al - N m - 1} e^{- \be t} g(t) dt < \infty \text{.} \label{eq:assumption2} 
\end{align}
Then, by Lemma \ref{lem:propriety}, the posterior under the prior $\p \sim \pi _{\al , \be , g, - m, \j ^{(m)}} ( \p )$ is proper, $K( \al , \be , g, r - m, \z + m \j ^{(N)} ) < \infty $ for all $\z \in {\mathbb{N} _0}^N$, and $K( \al + 1, \be , g, r - m, \z + m \j ^{(N)} ) < \infty $ for all $\z \in {\mathbb{N} _0}^N \setminus \{ \bm{0} ^{(N)} \} $. 

Define the function $\de ^{( \al , \be , g)} \colon {\mathbb{N} _0}^N \to (0, \infty ]$ by 
\begin{align}
\de ^{( \al , \be , g)} ( \z ) = {K( \al + 1, \be , g, r - m, \z + m \j ^{(N)} ) \over K( \al , \be , g, r - m, \z + m \j ^{(N)} )} \text{,} \quad \z \in {\mathbb{N} _0}^N \text{.} \non 
\end{align}
Let 
\begin{align}
\ph _{i, \nu }^{( \al , \be , g)} &= \begin{cases} \displaystyle {X_{i, \nu } \over r + X_{\cdot , \nu } - 1 + \de ^{( \al , \be , g)} ( \X _{\cdot } )} & \text{if $X_{i, \nu } \ge 1$} \\ \displaystyle 0 & \text{if $X_{i, \nu } = 0$} \end{cases} \non \\
&= {X_{i, \nu } \over r + X_{\cdot , \nu } - 1 + \de ^{( \al , \be , g)} ( \X _{\cdot } )} \label{eq:HB} 
\end{align}
for $i = 1, \dots , m$ and $\nu = 1, \dots , N$ and let $\hat{\p } ^{( \al , \be , g)} = ( \ph _{i, \nu }^{( \al , \be , g)} )_{1 \le i \le m, \, 1 \le \nu \le N}$. 
Then $\hat{\p } ^{( \al , \be , g)}$ is our hierarchical Bayes estimator. 

\begin{lem}
\label{lem:HB} 
Suppose that (\ref{eq:assumption1}) or (\ref{eq:assumption2}) holds. 
Then the shrinkage estimator $\hat{\p } ^{( \al , \be , g)}$ is a Bayes solution with respect to the prior (\ref{eq:shrinkage_prior}) with $a_0 = - m$ and $\a = \j ^{(m)}$ under the loss (\ref{eq:loss_SS}) for every $n = 1, \dots , N$. 
\end{lem}

\noindent
The term $\de ^{( \al , \be , g)} ( \X _{\cdot } )$ is at once expressed in closed form and symmetric in $X_{\cdot , 1} , \dots , X_{\cdot , N}$. 
Deriving such terms will be less straightforward in the case of empirical Bayes estimation except for those that are dependent only on $X_{\cdot , \cdot }$. 

Let $\e _{\nu }^{(N)}$ denote the $\nu $th unit vector in $\mathbb{R} ^N$, namely the $\nu $th column of the $N \times N$ identity matrix, for $\nu = 1, \dots , N$. 
The function $\de ^{( \al , \be , g)}$ satisfies the following properties. 
\begin{prp}
\label{prp:HB} 
Let $\z = ( z_1 , \dots , z_N )' \in {\mathbb{N} _0}^N$ and suppose that (\ref{eq:assumption1}) or (\ref{eq:assumption2}) holds. 
\begin{enumerate}
\item
We have $0 < \de ^{( \al , \be , g)} ( \z ) \le \infty $. 
Furthermore, $\de ^{( \al , \be , g)} ( \z ) = \infty $ only if $\z = \bm{0} ^{(N)}$. 
\item
Let $\nu = 1, \dots , N$. 
Then $\de ^{( \al , \be , g)} ( \z ) \ge \de ^{( \al , \be , g)} ( \z + \e _{\nu }^{(N)} )$. 
\item
Let $\nu = 1, \dots , N$. 
Then $\lim_{\mathbb{N} \ni k \to \infty } \de ^{( \al , \be , g)} ( \z + k \e _{\nu }^{(N)} ) = 0$. 
\item
Suppose that $r > m$, that $\lim_{t \to 0} g(t) = g(0) \in (0, \infty )$, and that $\al + 1 < N$. 
Then $\lim_{\mathbb{N} \setminus \{ 1 \} \ni k \to \infty } [ \de ^{( \al , \be , g)} ( \z + k \j ^{(N)} ) / \{ ( \al / N) / \log k \} ] = 1$. 
\end{enumerate}
\end{prp}

\noindent
Properties (iii) and (iv) above are in contrast to the fact that $\lim_{z \to \infty } \de ^{\rm{EB}} (z) = 1 + m > 0$.

The following theorem provides a sufficient condition for $\hat{\p } ^{( \al , \be , g)}$ to dominate $\hat{\p } ^{\rm{U}}$. 

\begin{thm}
\label{thm:HB} 
Let $n = 1, \dots , N$. 
Assume that (\ref{eq:assumption1}) or (\ref{eq:assumption2}) holds. 
Assume that $g$ is nonincreasing. 
Suppose further that %
\begin{align}
\al + 1 \le \min \{ n (m - 2), n m / 2 + \be r \} \text{.} \label{eq:assumption_HB} 
\end{align}
Then $\hat{\p } ^{( \al , \be , g)}$ is a hierarchical Bayes estimator dominating the UMVU estimator $\hat{\p } ^{\rm{U}}$ under the loss $L_n ( \d , \p )$ given by (\ref{eq:loss_SS}). 
\end{thm}

\noindent
There exist $\al > 0$ and $\be \ge 0$ satisfying assumption (\ref{eq:assumption_HB}) if and only if $n (m - 2) > 1$. 
When $r = m$ and $g$ is nonincreasing, the condition $\int_{0}^{1} t^{\al - N - 1} e^{- \be t} g(t) dt < \infty $ becomes $\al > N$. 
Even if $r = m$, the conditions of Theorem \ref{thm:HB} can be satisfied when $m$ is sufficiently large. 

In the remainder of this subsection, we consider the problem of estimating $\p $ under the loss (\ref{eq:loss_KL}) in order to show some robustness of our prior. 
Since the risk function of the UMVU estimator $\hat{\p }^{\rm{U}}$ is not defined under the loss (\ref{eq:loss_KL}), we first derive the Jeffreys prior. 

\begin{lem}
\label{lem:Jeff} 
The Dirichlet prior (\ref{eq:Dir}) with $a_0 = (1 - m) / 2$ and $\a = \j ^{(m)} / 2$ is the Jeffreys prior. 
\end{lem}

Next we show that under the loss (\ref{eq:loss_KL}), Bayes estimators are obtained as posterior means of $\p $. 

\begin{lem}
\label{lem:KL} 
Let $\p \sim \pi ( \p )$ be a strictly positive prior density and assume that the posterior is proper, that is, that $\int_{{D_m}^N} \big\{ \prod_{\nu = 1}^{N} {\rm{NM}}_m ( \x _{\nu } | r, \p _{\nu } ) \big\} \pi ( \p ) d\p < \infty $ for all $\x _1 , \dots , \x _N \in {\mathbb{N} _0}^m$. 
Then the posterior mean of $\p $ is a Bayes solution under the loss (\ref{eq:loss_KL}) for every $n = 1, \dots , N$. 
\end{lem}

The posterior under the Dirichlet prior (\ref{eq:Dir}) is proper if and only if $r + a_0 > 0$, in which case the posterior mean of $\p $ is 
\begin{align}
\hat{\p } ^{( a_0 , \a )} &= ( \ph _{i, \nu }^{( a_0 , \a )} )_{1 \le i \le m, \, 1 \le \nu \le N} \non \\
&= \Big( {X_{i, \nu } + a_i \over r + a_0 + X_{\cdot , \nu } + a_{\cdot }} \Big) _{1 \le i \le m, \, 1 \le \nu \le N} \text{.} \non 
\end{align}
The posterior under the hierarchical prior (\ref{eq:shrinkage_prior}) is proper if and only if the condition of part (ii) of Lemma \ref{lem:propriety} is satisfied. 
In this case, the posterior mean of $\p $ is 
\begin{align}
\hat{\p } ^{( \al , \be , g, a_0 , \a )} &= ( \ph _{i, \nu }^{( \al , \be , g, a_0 , \a )} )_{1 \le i \le m, \, 1 \le \nu \le N} \non \\
&= \Big( {X_{i, \nu } + a_i \over r + a_0 + X_{\cdot , \nu } + a_{\cdot } + \de _{\nu }^{( \al , \be , g, a_0 , \a )} ( \X _{\cdot } )} \Big) _{1 \le i \le m, \, 1 \le \nu \le N} \text{,} \non 
\end{align}
where $\de _{\nu }^{( \al , \be , g, a_0 , \a )} \colon {\mathbb{N} _0}^{N} \to (0, \infty )$ is the function defined by 
\begin{align}
\de _{\nu }^{( \al , \be , g, a_0 , \a )} ( \z ) &= {K( \al + 1, \be , g, r + a_0 , \z + a_{\cdot } \j ^{(N)} + \e _{\nu }^{(N)} ) \over K( \al , \be , g, r + a_0 , \z + a_{\cdot } \j ^{(N)} + \e _{\nu }^{(N)} )} \text{,} \quad \z \in {\mathbb{N} _0}^{N} \text{,} \non 
\end{align}
for $\nu = 1, \dots , N$. 
Some properties of the functions $\de _{\nu }^{( \al , \be , g, a_0 , \a )}$, $\nu = 1, \dots , N$, are given in the following proposition, which corresponds to Proposition \ref{prp:HB}. 

\begin{prp}
\label{prp:KL} 
Let $\z = ( z_1 , \dots , z_N )' \in {\mathbb{N} _0}^N$ and $\nu = 1, \dots , N$. 
Suppose that the condition of part (ii) of Lemma \ref{lem:propriety} is satisfied. 
\begin{enumerate}
\item
We have $0 < \de _{\nu }^{( \al , \be , g, a_0 , \a )} ( \z ) < \infty $. 
\item
Let ${\nu }' = 1, \dots , N$. 
Then $\de _{\nu }^{( \al , \be , g, a_0 , \a )} ( \z ) \ge \de _{\nu }^{( \al , \be , g, a_0 , \a )} ( \z + \e _{{\nu }'}^{(N)} )$. 
\item
Let ${\nu }' = 1, \dots , N$. 
Then $\lim_{\mathbb{N} \ni k \to \infty } \de _{\nu }^{( \al , \be , g, a_0 , \a )} ( \z + k \e _{{\nu }'}^{(N)} ) = 0$. 
\item
Suppose that $r + a_0 > 0$, that $\lim_{t \to 0} g(t) = g(0) \in (0, \infty )$, and that $\al + 1 < N$. 
Then $\lim_{\mathbb{N} \setminus \{ 1 \} \ni k \to \infty } \de _{\nu }^{( \al , \be , g, a_0 , \a )} ( \z + k \j ^{(N)} ) / \{ ( \al / N) / \log k \} ] = 1$. 
\end{enumerate}
\end{prp}

Theorem \ref{thm:KL} provides a sufficient condition for $\hat{\p } ^{( \al , \be , g, a_0 , \a )}$ to dominate $\hat{\p } ^{( a_0 , \a )}$ under the loss (\ref{eq:loss_KL}). 

\begin{thm}
\label{thm:KL} 
Let $n = 1, \dots , N$. 
Assume that the condition of part (ii) of Lemma \ref{lem:propriety} is satisfied. 
Assume that $g$ is nonincreasing. 
Suppose further that $a_0 + a_{\cdot } + 1 \ge 0$ and that 
\begin{align}
\al + 1 \le n (- a_0 - 2 ) \text{.} \label{eq:assumption_KL} 
\end{align}
Then $\hat{\p } ^{( \al , \be , g, a_0 , \a )}$ dominates $\hat{\p } ^{( a_0 , \a )}$ under the loss $\tilde{L} _n ( \tilde{\d } , \p )$ given by (\ref{eq:loss_KL}). 
\end{thm}

In particular, we have the following result for the case of the Jeffreys prior. 
\begin{cor}
\label{cor:Jeff} 
Let $n = 1, \dots , N$. 
Assume that either 
\begin{itemize}
\item
$r > (m - 1)/ 2$ and $\int_{1}^{\infty } t^{\al - N a_{\cdot } - 1} e^{- \be t} g(t) dt < \infty $ or 
\item
$r = (m - 1)/ 2$, $\al > N$, and $\int_{1}^{\infty } t^{\al - N a_{\cdot } - 1} e^{- \be t} g(t) dt < \infty $. 
\end{itemize}
Assume that $g$ is nonincreasing. 
Suppose further that 
\begin{align}
\al + 1 \le n (m - 5) / 2 \text{.} \label{eq:assumption_Jeff} 
\end{align}
Then $\hat{\p } ^{( \al , \be , g, (1 - m) / 2 , \j ^{(m)} / 2)}$ dominates $\hat{\p } ^{((1 - m) / 2 , \j ^{(m)} / 2)}$ under the loss $\tilde{L} _n ( \tilde{\d } , \p )$ given by (\ref{eq:loss_KL}). 
\end{cor}

\subsection{Posterior computation}
\label{subsec:hb_computation} 
In order to approximate the integral 
\begin{align}
K( \al , \be , g, \xi _0 , \bxi ) &= \int_{0}^{\infty } t^{\al - 1} e^{- \be t} g(t) \Big\{ \prod_{\nu = 1}^{N} {\Ga (t + \xi _0 ) \over \Ga (t + \xi _0 + \xi _{\nu } )} \Big\} dt \text{,} \non 
\end{align}
we could in principle use i.i.d. gamma variables (when $\be > 0$) or rewrite the integral as 
\begin{align}
K( \al , \be , g, \xi _0 , \bxi ) &= \int_{0}^{1} \Big( {\om \over 1 - \om } \Big) ^{\al - 1} e^{- \be \om / (1 - \om )} g \Big( {\om \over 1 - \om } \Big) \Big\{ \prod_{\nu = 1}^{N} {\Ga ( \om / (1 - \om ) + \xi _0 ) \over \Ga ( \om / (1 - \om ) + \xi _0 + \xi _{\nu } )} \Big\} {1 \over (1 - \om )^2} d\om \non 
\end{align}
and use i.i.d. uniform variables, for example. 
However, this can be numerically unstable because of the gamma function in the integrand. 
If $\bxi \in {\mathbb{N} _0}^N$, the problem would be alleviated to some extent by using the relation 
\begin{align}
\prod_{\nu = 1}^{N} {\Ga (t + \xi _0 ) \over \Ga (t + \xi _0 + \xi _{\nu } )} = \prod_{\nu = 1}^{N} {1 \over (t + \xi _0 ) \dotsm (t + \xi _0 + \xi _{\nu } - 1)} \non 
\end{align}
for all $t \in (0, \infty )$. 

When $g = g_1$, a more convenient way to compute the hierarchical Bayes estimators in the previous subsection is to use MCMC samples since they are functions of posterior expectations. 
In order to describe a Gibbs sampler, we introduce a fully conjugate prior. 
For $\al > 0$, $\be \ge 0$, $a_0 \in \mathbb{R}$, and $( \a _1 , \dots , \a _N ) = ( ( a_{1, 1} , \dots , a_{m, 1} )' , \dots , ( a_{1, N} , \dots , a_{m, N} )' ) \in (0, \infty )^{m \times N}$, let 
\begin{align}
\pi ( \p , t | \al , \be , a_0 , \a _1 , \dots , \a _N ) &= t^{\al - 1} e^{- \be t} \prod_{\nu = 1}^{N} \Big( {p_{0, \nu }}^{t + a_{0} - 1} \prod_{i = 1}^{m} {p_{i, \nu }}^{a_{i, \nu } - 1} \Big)  \label{eq:joint_prior} 
\end{align}
and 
\begin{align}
\pi ( \p | \al , \be , a_0 , \a _1 , \dots , \a _N ) &= \Ga ( \al ) \Big\{ \prod_{\nu = 1}^{N} \Big( {p_{0, \nu }}^{a_{0} - 1} \prod_{i = 1}^{m} {p_{i, \nu }}^{a_{i, \nu } - 1} \Big) \Big\} / \Big( \be + \sum_{\nu = 1}^{N} \log {1 \over p_{0, \nu }} \Big) ^{\al } \text{,} \label{eq:marginal_prior} 
\end{align}
where $t \in (0, \infty )$ and where $\p = %
( ( p_{1, 1} , \dots , p_{m, 1} )' , \dots , ( p_{1, N} , \dots , p_{m, N} )' ) \in {D_m}^N$ and $p_{0, \nu } %
= 1 - \sum_{i = 1}^{m} p_{i, \nu }$ for $\nu = 1, \dots , N$. 
When $\a _1 = \dots = \a _N = \a $, the prior (\ref{eq:marginal_prior}) becomes the original prior (\ref{eq:shrinkage_prior_simple}). 

Some basic properties of the priors (\ref{eq:joint_prior}) and (\ref{eq:marginal_prior}) are summarized in the following proposition. 
Let $a_{\cdot , \nu } = \sum_{i = 1}^{m} a_{i, \nu }$ for $\nu = 1, \dots , N$ and let $a_{\cdot , \cdot } = \sum_{\nu = 1}^{N} a_{\cdot , \nu }$. 

\begin{prp}
\label{prp:MCMC} 
The priors (\ref{eq:joint_prior}) and (\ref{eq:marginal_prior}) satisfy the following properties: 
\begin{enumerate}
\item
The following are equivalent: 
\begin{itemize}
\item
$\int_{{D_m}^N \times (0, \infty )} \pi ( \p , t | \al , \be , a_0 , \a _1 , \dots , \a _N ) d( \p , t) < \infty $. 
\item
$\int_{{D_m}^N} \pi ( \p | \al , \be , a_0 , \a _1 , \dots , \a _N ) d\p < \infty $. 
\item
$\min \{ \max \{ a_0 , \al - N \} , \max \{ a_{\cdot , \cdot } - \al , \be \} \} > 0$. 
\end{itemize}
\item
If $\p \sim \pi ( \p | \al , \be , a_0 , \a _1 , \dots , \a _N )$ and $( \x _1 , \dots , \x _N ) | \p \sim \prod_{\nu = 1}^{N} {\rm{NM}}_m ( \x _{\nu } | r, \p _{\nu } )$, then %
\begin{align}
\p | ( \x _1 , \dots , \x _N ) \sim \pi ( \p | \al , \be , r + a_0 , \x _1 + \a _1 , \dots , \x _N + \a _N ) \text{.} \non 
\end{align}
\item
If $( \p , t) \sim \pi ( \p , t | \al , \be , a_0 , \a _1 , \dots , \a _N )$, then $\p \sim \pi ( \p | \al , \be , a_0 , \a _1 , \dots , \a _N )$. 
\item
If $( \p , t) \sim \pi ( \p , t | \al , \be , a_0 , \a _1 , \dots , \a _N )$, then %
\begin{align}
\p | t &\sim \prod_{\nu = 1}^{N} {\rm{Dir}}_m ( \p _{\nu } | t + a_0 , \a _{\nu } ) \text{,} \non \\
t | \p &\sim {\rm{Ga}} \Big( t \Big| \al , \be + \sum_{\nu = 1}^{N} \log {1 \over p_{0, \nu }} \Big) \text{.} \non 
\end{align}
\end{enumerate}
\end{prp}

Part (ii) of Proposition \ref{prp:MCMC} shows that the prior (\ref{eq:marginal_prior}) is conjugate. 
Furthermore, part (iii) of the proposition shows that in order to generate samples of $\p $ from the prior (\ref{eq:marginal_prior}), it is sufficient to sample from the joint prior (\ref{eq:joint_prior}). 
Therefore, we describe a Gibbs sampler for (\ref{eq:joint_prior}) based on part (iv) of the proposition. 
In order to generate MCMC samples corresponding to (\ref{eq:joint_prior}) when it is proper, given a current sample of $( \p , t)$, denoted by $( \tilde{\p } , \tilde{t} ) = ( ( \tilde{p} _{i, \nu } )_{1 \le i \le m, \, 1 \le \nu \le N} , \tilde{t} )$, we generate a new sample as follows: 
\begin{itemize}
\item
sample $t^{*} \sim {\rm{Ga}} \big( t \big| \al , \be + \sum_{\nu = 1}^{N} \log \big\{ 1 / \big( 1 - \sum_{i = 1}^{m} \tilde{p} _{i, \nu } \big) \big\} \big) $; 
\item
sample $\p ^{*} \sim \prod_{\nu = 1}^{N} {\rm{Dir}}_m ( \p _{\nu } | t^{*} + a_0 , \a _{\nu } )$. 
\end{itemize}
Then samples of $\p $ can be used to approximate expectations of functions of $\p \sim \pi ( \p | \al , \be , a_0 , \a _1 , \dots , \a _N )$. 
Also, samples of $t$ may be used to approximate $\de ^{( \al , \be , g)}$ and $\de _{\nu }^{( \al , \be , g, a_0 , \a )}$, $\nu = 1, \dots , N$, even if $g \neq g_1$.

\section{Simulation Study}
\label{sec:sim}
In this section, we investigate through simulation the numerical performance of the risk functions of the Bayes estimators given in the previous two sections under the standardized squared error loss given by (\ref{eq:loss_SS}) with $n = N$. 
The estimators which we compare are the following four: 

\medskip
U: the UMVU estimator $\hat{\p } ^{\rm{U}}$ given by (\ref{eq:UMVU}), 

\medskip
EB0: the empirical Bayes estimator which estimates $\p _1 , \dots , \p _N$ independently based on $\X _1 , \dots , \X _N$, respectively, namely $\hat{\p } ^{\rm{EB0}} = ( X_{i, \nu } / (r + X_{\cdot , \nu } + m + m r / X_{\cdot , \nu } ) )_{1 \le i \le m, \, 1 \le \nu \le N}$, 

\medskip
EB: the empirical Bayes estimator $\hat{\p } ^{\rm{EB}}$ given by (\ref{eq:EB}), 

\medskip
HB: the hierarchical Bayes estimator $\hat{\p } ^{\rm{HB}} = \hat{\p } ^{( \al , 1, g_1 )}$ given by (\ref{eq:HB}) with $( \be , g) = (1, g_1 )$. 

\medskip
We consider the following cases: 
\begin{itemize}
\item[(i)]
We set $(r, m, N) = (8, 7, 3)$, $\al = 14$, and $\p = \p ^{(1)} (1), \p ^{(1)} (2), \p ^{(1)} (3)$, where 
\begin{align}
\p ^{(1)} (1) &= ((1, 1, 1, 1, 1, 1, 1 )' / 8, (1, 1, 1, 1, 1, 1, 1 )' / 8, (1, 1, 1, 1, 1, 1, 1 )' / 8) \text{,} \non \\
\p ^{(1)} (2) &= ((1, 1, 1, 1, 2, 2, 2 )' / 12, (1, 1, 1, 1, 1, 1, 1 )' / 8, (1, 1, 1, 1, 2, 2, 2 )' / 12) \text{,} \non \\
\p ^{(1)} (3) &= ((1, 1, 1, 1, 2, 2, 2 )' / 12, (1, 1, 1, 1, 1, 1, 1 )' / 8, (2, 2, 2, 2, 1, 1, 1 )' / 12) \text{.} \non 
\end{align}
\item[(ii)]
We set $(r, m, N) = (4, 3, 7)$, $\al = 6$, and $\p = \p ^{(2)} (1), \p ^{(2)} (2), \p ^{(2)} (3)$, where 
\begin{align}
\p ^{(2)} (1) &= \begin{pmatrix} \begin{pmatrix} 1 \\ 1 \\ 1 \end{pmatrix} / 4 & \begin{pmatrix} 1 \\ 1 \\ 1 \end{pmatrix} / 4 & \begin{pmatrix} 1 \\ 1 \\ 1 \end{pmatrix} / 4 & \begin{pmatrix} 1 \\ 1 \\ 1 \end{pmatrix} / 4 & \begin{pmatrix} 1 \\ 1 \\ 1 \end{pmatrix} / 4 &\begin{pmatrix} 1 \\ 1 \\ 1 \end{pmatrix} / 4 & \begin{pmatrix} 1 \\ 1 \\ 1 \end{pmatrix} / 4 \end{pmatrix} \text{,} \non \\
\p ^{(2)} (2) &= \begin{pmatrix} \begin{pmatrix} 1 \\ 1 \\ 2 \end{pmatrix} / 6 & \begin{pmatrix} 1 \\ 1 \\ 2 \end{pmatrix} / 6 & \begin{pmatrix} 1 \\ 1 \\ 1 \end{pmatrix} / 4 & \begin{pmatrix} 1 \\ 1 \\ 1 \end{pmatrix} / 4 & \begin{pmatrix} 1 \\ 1 \\ 1 \end{pmatrix} / 4 & \begin{pmatrix} 1 \\ 1 \\ 2 \end{pmatrix} / 6 & \begin{pmatrix} 1 \\ 1 \\ 2 \end{pmatrix} / 6 \end{pmatrix} \text{,} \non \\
\p ^{(2)} (3) &= \begin{pmatrix} \begin{pmatrix} 1 \\ 1 \\ 2 \end{pmatrix} / 6 & \begin{pmatrix} 1 \\ 1 \\ 2 \end{pmatrix} / 6 & \begin{pmatrix} 1 \\ 1 \\ 1 \end{pmatrix} / 4 & \begin{pmatrix} 1 \\ 1 \\ 1 \end{pmatrix} / 4 & \begin{pmatrix} 1 \\ 1 \\ 1 \end{pmatrix} / 4 & \begin{pmatrix} 2 \\ 2 \\ 1 \end{pmatrix} / 6 & \begin{pmatrix} 2 \\ 2 \\ 1 \end{pmatrix} / 6 \end{pmatrix} \text{.} \non 
\end{align}
\item[(iii)]
We set $(r, m, N) = (2, 1, 7)$, $\al = 6$, and $\p = \p ^{(3)} (1), \p ^{(3)} (2), \p ^{(3)} (3)$, where 
\begin{align}
\p ^{(3)} (1) &= (1 / 2, 1 / 2, 1 / 2, 1 / 2, 1 / 2, 1 / 2, 1 / 2) \text{,} \non \\
\p ^{(3)} (2) &= (1 / 3, 1 / 3, 1 / 2, 1 / 2, 1 / 2, 1 / 3, 1 / 3) \text{,} \non \\
\p ^{(3)} (3) &= (1 / 3, 1 / 3, 1 / 2, 1 / 2, 1 / 2, 2 / 3, 2 / 3) \text{.} \non 
\end{align}
\end{itemize}
Case (i) is a case where $m > N$ while case (ii) is where $m < N$. 
Case (iii) corresponds to the negative binomial distribution. 
Table \ref{table:condition} summarizes whether the sufficient conditions for dominance in Sections \ref{sec:eb} and \ref{sec:hb} are applicable.

\small
\begin{table}[!thb]
\caption{Whether the conditions for dominance are satisfied or not. (When one of the conditions is satisficed, $+$ is marked, and $-$ is marked otherwise.) }
\begin{center}
$
{\renewcommand\arraystretch{1.1}\small
\begin{array}{c@{\hspace{5mm}}
              r@{\hspace{2mm}}
              r@{\hspace{2mm}}
              r@{\hspace{2mm}}
              r@{\hspace{2mm}}
              r@{\hspace{2mm}}
              r@{\hspace{2mm}}
              r@{\hspace{2mm}}
              r
             }
\text{Case} & \text{EB0} & \text{EB} & \text{HB}\\

\hline

\text{(i)}
&$$
{}+{}
$$
&$$
{}+{}
$$
&$$
{}+{}
$$

\\

\text{(ii)}
&$$
{}-{}
$$
&$$
{}-{}
$$
&$$
{}+{}
$$

\\

\text{(iii)}
&$$
{}-{}
$$
&$$
{}-{}
$$
&$$
{}-{}
$$

\\

\end{array}
}
$
\end{center}
\label{table:condition}
\end{table}
\normalsize

For each estimator $\hat{\p }$, we obtain approximated values of the risk function $E[ L_N ( \hat{\p } , \p ) ] $ by simulation with $1,000$ replications. 
The hierarchical Bayes estimator $\hat{\p } ^{\rm{HB}}$ was computed based on the Gibbs sampler described in Section \ref{subsec:hb_computation} by generating $50,000$ posterior samples after discarding the first $50,000$ samples. 
The percentage relative improvement in average loss (PRIAL) of an estimator $\hat{\p }$ over $\hat{\p } ^{\rm{U}}$ is defined by 
\begin{align}
{\rm PRIAL} = 100 \{ E[ L_N ( \hat{\p } ^{\rm{U}} , \p ) ] - E[ L_N ( \hat{\p } , \p ) ] \} / E[ L_N ( \hat{\p } ^{\rm{U}} , \p ) ] \text{.} \non 
\end{align}

For case (i), Table \ref{table:risk_1} reports values of the risks of the estimators with values of PRIAL given in parentheses. 
In all cases, the risk values of $\hat{\p } ^{\rm{EB0}}$ are smaller than those of $\hat{\p } ^{\rm{HB}}$, and the risk values of $\hat{\p } ^{\rm{EB}}$ are still smaller. 
These three estimators have the largest values of PRIAL when $\p = \p ^{(1)} (2)$. 
Also, it can be seen that in the balanced case of $\p = \p ^{(1)} (1)$, the risk values of the three estimators are smaller than those of $\hat{\p } ^{\rm{U}}$ even when the loss is (\ref{eq:loss_SS}) with $n = 1$.

\small
\begin{table}[!thb]
\caption{Risks of the estimators U, EB0, EB, and HB for case (i). (Values of PRIAL of EB0, EB, and HB are given in parentheses) }
\begin{center}
$
{\renewcommand\arraystretch{1.1}\small
\begin{array}{c@{\hspace{5mm}}
              r@{\hspace{2mm}}
              r@{\hspace{2mm}}
              r@{\hspace{2mm}}
              r@{\hspace{2mm}}
              r@{\hspace{2mm}}
              r@{\hspace{2mm}}
              r@{\hspace{2mm}}
              r
             }
\text{$\p $} & \text{U} & \text{EB0} & \text{EB} & \text{HB} \\

\hline

\text{$\p ^{(1)} (1)$}
&0.32 	

&0.30  \, ( 7.18 	
 ) 
& 	0.29 	\, ( 7.91 	  ) 
&0.31  \, ( 3.70 )

\\
\text{$\p ^{(1)} (2)$}
&0.39 	

& 0.35 	 \, (  9.65 	
) 
&0.35 	\, ( 10.61  ) 
&0.37  \, ( 	5.17  )

\\
\text{$\p ^{(1)} (3)$}
&0.32 	

&0.29 	 \, (  9.10 
) 
&0.29   \, ( 	9.92 ) 
&	0.31 \, ( 	4.15   )

\\

\end{array}
}
$
\end{center}
\label{table:risk_1}
\end{table}
\normalsize

For Case %
(ii), Table \ref{table:risk_2} reports values of  the risks and PRIAL. 
Although the empirical Bayes estimators do not satisfy the condition of Corollary \ref{cor:EB}, 
$\hat{\p } ^{\rm{EB0}}$ is competitive with $\hat{\p } ^{\rm{HB}}$ and $\hat{\p } ^{\rm{EB}}$ is superior to $\hat{\p } ^{\rm{HB}}$. 

\small
\begin{table}[!thb]
\caption{Risks of the estimators U, EB0, EB, and HB for case (ii). (Values of PRIAL of EB0, EB, and HB are given in parentheses) }
\begin{center}
$
{\renewcommand\arraystretch{1.1}\small
\begin{array}{c@{\hspace{5mm}}
              r@{\hspace{2mm}}
              r@{\hspace{2mm}}
              r@{\hspace{2mm}}
              r@{\hspace{2mm}}
              r@{\hspace{2mm}}
              r@{\hspace{2mm}}
              r@{\hspace{2mm}}
              r
             }
\text{$\p $} & \text{U} & \text{EB0} & \text{EB} & \text{HB} \\

\hline

\text{$\p ^{(2)} (1)$}
&1.21 	

&1.16  \, (4.25 
  ) 
& 	1.10 	 \, (  	9.27 	 ) 
&1.14 \, ( 6.16 ) 

\\
\text{$\p ^{(2)} (2)$}
&1.44 	

&1.30 	 \, ( 9.90 
 ) 
&1.23 	\, ( 	14.82 	  ) 
&1.32   \, ( 8.55 ) 

\\
\text{$\p ^{(2)} (3)$}
&1.22 	

& 1.15 	\, ( 6.21 	
 ) 
&1.08 	\, ( 11.58 	) 
&1.14   \, (  6.87  )

\\

\end{array}
}
$
\end{center}
\label{table:risk_2}
\end{table}
\normalsize

Finally, Table \ref{table:risk_3} reports values of  the risks and PRIAL for case (iii). 
The estimators $\hat{\p } ^{\rm{EB}}$ and $\hat{\p } ^{\rm{HB}}$ do not satisfy the conditions for dominance but their risk values are smaller than those of $\hat{\p } ^{\rm{U}}$. 
In particular, $\hat{\p } ^{\rm{HB}}$ has large values of PRIAL.

\small
\begin{table}[!thb]
\caption{Risks of the estimators U, EB0, EB, and HB for case (iii). (Values of PRIAL of EB0, EB, and HB are given in parentheses) }
\begin{center}
$
{\renewcommand\arraystretch{1.1}\small
\begin{array}{c@{\hspace{5mm}}
              r@{\hspace{2mm}}
              r@{\hspace{2mm}}
              r@{\hspace{2mm}}
              r@{\hspace{2mm}}
              r@{\hspace{2mm}}
              r@{\hspace{2mm}}
              r@{\hspace{2mm}}
              r
             }
\text{$\p $} & \text{U} & \text{EB0} & \text{EB} & \text{HB} \\

\hline

\text{$\p ^{(3)} (1)$}
&1.34 	

&1.38 	\, ( -3.35 	
 ) 
&1.33   \, ( 0.75  ) 
&	1.00  \, (	24.99   )

\\
\text{$\p ^{(3)} (2)$}
&1.72 	

&1.32 	  \, (  23.43 	
) 
&1.30  \, ( 24.39  ) 
&	1.11 \, (  	35.47 )

\\
\text{$\p ^{(3)} (3)$}
&1.36 	

&1.28 	  \, ( 5.78 
) 
&1.23  \, ( 	9.62   ) 
&	0.97 \, (  	28.42)

\\

\end{array}
}
$
\end{center}
\label{table:risk_3}
\end{table}
\normalsize

\section{Appendix}
Here we give proofs. 
Let $\bm{0} ^{(m)} = (0, \dots , 0)' \in \mathbb{R} ^m$ and $\bm{0} ^{(m, N)} = \bm{0} ^{(m)} {\bm{0} ^{(N)}}' \in \mathbb{R} ^{m \times N}$. 
Let $\e _{i}^{(m)}$ be the $i$th unit vector in $\mathbb{R} ^m$, namely the $i$th column of the $m \times m$ identity matrix, for $i = 1, \dots , m$. 
Let $\e _{i, \nu }^{(m, N)} = \e _{i}^{(m)} ( \e _{\nu }^{(N)} )' \in \mathbb{R} ^{m \times N}$ for $i = 1, \dots , m$ and $\nu = 1, \dots , N$. 
Further let $\de _{i, j}^{(m)} = {\e _{i}^{(m)}}' \e _{j}^{(m)}$ for $i, j = 1, \dots , m$ and let $\de _{\nu, {\nu }'}^{(N)} = {\e _{\nu }^{(N)}}' \e _{{\nu }'}^{(N)}$ for $\nu , {\nu }' = 1, \dots , N$. 
For $\v = ( v_{i, \nu } )_{1 \le i \le m, \, 1 \le \nu \le N} \in \mathbb{R} ^{m \times N}$ and $\tilde{\v } = ( \tilde{v} _{i, \nu } )_{1 \le i \le m, \, 1 \le \nu \le N} \in \mathbb{R} ^{m \times N}$, we write the inner product $\sum_{\nu = 1}^{N} \sum_{i = 1}^{m} v_{i, \nu } \tilde{v} _{i, \nu }$ as $\v \cdot \tilde{\v }$. 
The following result is due to Hudson (1978). 

\begin{lem}
\label{lem:hudson}
Let $h \colon {\mathbb{N} _0}^{m \times N} \to \mathbb{R}$ and suppose that either $h( \x ) \ge 0$ for all $\x \in {\mathbb{N} _0}^{m \times N}$ or $E[ |h( \X )| ] < \infty $. 
Then for all $i = 1, \dots , m$ and all $\nu = 1, \dots , N$, if $h( \x ) = 0$ for all $\x = ( x_{j, {\nu }'} )_{1 \le j \le m, \, 1 \le {\nu }' \le N} \in {\mathbb{N} _0}^{m \times N}$ such that $x_{i, \nu } = 0$, we have 
\begin{align}
E \Big[ {h( \X ) \over p_{i, \nu }} \Big] = E \Big[ {r + X_{\cdot , \nu } \over X_{i, \nu } + 1} h( \X + \e _{i, \nu }^{(m, N)} ) \Big] \text{.} \non 
\end{align}
\end{lem}

\noindent
{\bf Proof%
.} \ \ We have 
\begin{align}
E \Big[ {h( \X ) \over p_{i, \nu }} \Big] &= \sum_{\x = ( \x _1 , \dots , \x _N ) \in {\mathbb{N} _0}^{m \times N} , \, \x \cdot \e _{i, \nu }^{(m, N)} \neq 0} {h( \x ) \over p_{i, \nu }} \prod_{{\nu }' = 1}^{N} {\rm{NM}}_m ( \x _{{\nu }'} | r, \p _{{\nu }'} ) \non \\
&= \sum_{\x = ( \x _1 , \dots , \x _N ) \in {\mathbb{N} _0}^{m \times N}} {h( \x + \e _{i, \nu }^{(m, N)} ) \over p_{i, \nu }} {{\rm{NM}}_m ( \x _{\nu } + \e _{i}^{(m)} | r, \p _{\nu } ) \over {\rm{NM}}_m ( \x _{\nu } | r, \p _{\nu } )} \prod_{{\nu }' = 1}^{N} {\rm{NM}}_m ( \x _{{\nu }'} | r, \p _{{\nu }'} ) \non \\
&= E \Big[ {r + X_{\cdot , \nu } \over X_{i, \nu } + 1} h( \X + \e _{i, \nu }^{(m, N)} ) \Big] \text{,} \non 
\end{align}
which proves the desired result. 
\hfill$\Box$

\bigskip

\noindent
{\bf Proof of Theorem \ref{thm:shrinkage_estimator_general}.} \ \ Let $\De _{n}^{( \de )} = E[ L_n ( \hat{\p } ^{( \de )} , \p ) ] - E[ L_n ( \hat{\p } ^{\rm{U}} , \p ) ]$. 
For $\nu = 1, \dots , N$, let 
\begin{align}
\phi _{\nu }^{( \de )} ( \X ) %
&= {\de ( X_{\cdot , \cdot } ) \over r + X_{\cdot , \nu } - 1 + \de ( X_{\cdot , \cdot } )} \non 
\end{align}
so that 
\begin{align}
\ph _{i, \nu }^{( \de )} = \ph _{i, \nu }^{\rm{U}} - \ph _{i, \nu }^{\rm{U}} \phi _{\nu }^{( \de )} ( \X ) \non 
\end{align}
for every $i = 1, \dots , m$. 
Then, by Lemma \ref{lem:hudson}, we have 
\begin{align}
\De _{n}^{( \de )} &= E \Big[ \sum_{\nu = 1}^{n} \sum_{i = 1}^{m} \Big( {1 \over p_{i, \nu }} [ ( \ph _{i, \nu }^{\rm{U}} )^2 \{ \phi _{\nu }^{( \de )} ( \X ) \} ^2 - 2 ( \ph _{i, \nu }^{\rm{U}} )^2 \phi _{\nu }^{( \de )} ( \X )] + 2 \ph _{i, \nu }^{\rm{U}} \phi _{\nu }^{( \de )} ( \X ) \Big) \Big] \non \\
&= E \Big[ \sum_{\nu = 1}^{n} \sum_{i = 1}^{m} \Big[ {X_{i, \nu } + 1 \over r + X_{\cdot , \nu }} \{ \phi _{\nu }^{( \de )} ( \X + \e _{i, \nu }^{(m, N)} ) \} ^2 - 2 {X_{i, \nu } + 1 \over r + X_{\cdot , \nu }} \phi _{\nu }^{( \de )} ( \X + \e _{i, \nu }^{(m, N)} ) + 2 \ph _{i, \nu }^{\rm{U}} \phi _{\nu }^{( \de )} ( \X ) \Big] \Big] \non \\
&= E \Big[ \sum_{\nu = 1}^{n} \{ I_{1, \nu }^{( \de )} ( \X ) - 2 I_{2, \nu }^{( \de )} ( \X ) + 2 I_{3, \nu }^{( \de )} ( \X ) \} \Big] \text{,} \non 
\end{align}
where 
\begin{align}
I_{1, \nu }^{( \de )} ( \x ) &= {\sum_{i = 1}^{m} x_{i, \nu } + m \over r + \sum_{i = 1}^{m} x_{i, \nu }} \Big\{ {\de \big( \sum_{{\nu }' = 1}^{N} \sum_{i = 1}^{m} x_{i, {\nu }'} + 1 \big) \over r + \sum_{i = 1}^{m} x_{i, \nu } + \de \big( \sum_{{\nu }' = 1}^{N} \sum_{i = 1}^{m} x_{i, {\nu }'} + 1 \big) } \Big\} ^2 \text{,} \non \\
I_{2, \nu }^{( \de )} ( \x ) &= {\sum_{i = 1}^{m} x_{i, \nu } + m \over r + \sum_{i = 1}^{m} x_{i, \nu }} {\de \big( \sum_{{\nu }' = 1}^{N} \sum_{i = 1}^{m} x_{i, {\nu }'} + 1 \big) \over r + \sum_{i = 1}^{m} x_{i, \nu } + \de \big( \sum_{{\nu }' = 1}^{N} \sum_{i = 1}^{m} x_{i, {\nu }'} + 1 \big) } \text{,} \non \\
I_{3, \nu }^{( \de )} ( \x ) %
&= {\sum_{i = 1}^{m} x_{i, \nu } \over r + \sum_{i = 1}^{m} x_{i, \nu } - 1} {\de \big( \sum_{{\nu }' = 1}^{N} \sum_{i = 1}^{m} x_{i, {\nu }'} \big) \over r + \sum_{i = 1}^{m} x_{i, \nu } - 1 + \de \big( \sum_{{\nu }' = 1}^{N} \sum_{i = 1}^{m} x_{i, {\nu }'} \big) } \text{,} \non 
\end{align}
for $\x = %
( x_{i, {\nu }'} )_{1 \le i \le m, \, 1 \le {\nu }' \le N} \in {\mathbb{N} _0}^{m \times N}$ for each $\nu = 1, \dots , N$. 
Since $\sum_{\nu = 1}^{n} \{ I_{1, \nu }^{( \de )} ( \bm{0} ^{(m, N)} ) - 2 I_{2, \nu }^{( \de )} ( \bm{0} ^{(m, N)} ) + 2 I_{3, \nu }^{( \de )} ( \bm{0} ^{(m, N)} ) \} < 0$, it is sufficient to show that $\sum_{\nu = 1}^{n} \{ I_{1, \nu }^{( \de )} ( \x ) - 2 I_{2, \nu }^{( \de )} ( \x ) + 2 I_{3, \nu }^{( \de )} ( \x ) \} \le 0$ for all $\x \in {\mathbb{N} _0}^{m \times N} \setminus \{ \bm{0} ^{(m, N)} \} $. 

Fix $\x = %
( x_{i, \nu } )_{1 \le i \le m, \, 1 \le \nu \le N} \in {\mathbb{N} _0}^{m \times N} \setminus \{ \bm{0} ^{(m, N)} \} $. 
For notational simplicity, let $z _{\nu } = \sum_{i = 1}^{m} x_{i, \nu }$ for $\nu = 1, \dots , N$ and let $z = \sum_{\nu = 1}^{N} z_{\nu }$. 
Then for all $\nu = 1, \dots , N$ such that $z_{\nu } \neq 0$, since 
\begin{align}
\de (z) &\le {z + 1 \over z} \de (z + 1) \le {z_{\nu } + 1 \over z_{\nu }} \de (z + 1) \text{,} \non 
\end{align}
we have 
\begin{align}
I_{3, \nu }^{( \de )} ( \x ) &= {z_{\nu } \over r + z_{\nu } - 1} {\de (z) \over r + z_{\nu } - 1 + \de (z)} \non \\
&\le {z_{\nu } \over r + z_{\nu } - 1} {( z_{\nu } + 1) \de (z + 1) \over z_{\nu } (r + z_{\nu } - 1) + ( z_{\nu } + 1) \de (z + 1)} \non \\
&\le {z_{\nu } + 3 \over r + z_{\nu }} {\de (z + 1) \over r + z_{\nu } + \de (z + 1)} \text{,} \non 
\end{align}
where the second inequality follows from the assumption that $r \ge 5 / 2$. 
Therefore, 
\begin{align}
\sum_{\nu = 1}^{n} \{ I_{1, \nu }^{( \de )} ( \x ) - 2 I_{2, \nu }^{( \de )} ( \x ) + 2 I_{3, \nu }^{( \de )} ( \x ) \} \le I^{( \de )} ( \x ) \text{,} \non 
\end{align}
where 
\begin{align}
I^{( \de )} ( \x ) &= \sum_{\nu  = 1}^{n} \Big( {1 \over r + z_{\nu }} {\de (z + 1) \over \{ r + z_{\nu } + \de (z + 1) \} ^2} \non \\
&\quad \times [ z_{\nu } \{ \de (z + 1) - 2 (m - 3) \} - (m - 6) \de (z + 1) - 2 (m - 3) r] \Big) \text{.} \non 
\end{align}
Suppose first that $\de (z + 1) \le 2 (m - 3)$. 
Then $I^{( \de )} ( \x ) \le 0$ by assumption since $z + 1 \ge 2$. 
On the other hand, if $\de (z + 1) > 2 (m - 3)$, then, by the covariance inequality, 
\begin{align}
I^{( \de )} ( \x ) &\le {1 \over n} \Big[ \sum_{\nu = 1}^{n} {1 \over r + z_{\nu }} {\de (z + 1) \over \{ r + z_{\nu } + \de (z + 1) \} ^2} \Big] \non \\
&\quad \times \Big[ \Big( \sum_{\nu = 1}^{n} z_{\nu } \Big) \{ \de (z + 1) - 2 (m - 3) \} - n \{ (m - 6) \de (z + 1) + 2 (m - 3) r \} \Big] \non \\
&\le {1 \over n} \Big[ \sum_{\nu = 1}^{n} {1 \over r + z_{\nu }} {\de (z + 1) \over \{ r + z_{\nu } + \de (z + 1) \} ^2} \Big] \non \\
&\quad \times [z \{ \de (z + 1) - 2 (m - 3) \} - n \{ (m - 6) \de (z + 1) + 2 (m - 3) r \} ] \text{.} \non 
\end{align}
The right-hand side of the above inequality is nonpositive by assumption. 
This completes the proof. 
\hfill$\Box$

\begin{remark}
\label{rem:neccesary_condition} 
In the above proof, we have shown that $I_{n}^{( \de )} ( \x ) = \sum_{\nu = 1}^{n} \{ I_{1, \nu }^{( \de )} ( \x ) - 2 I_{2, \nu }^{( \de )} ( \x ) + 2 I_{3, \nu }^{( \de )} ( \x ) \} \le 0$ for all $\x \in {\mathbb{N} _0}^{m \times N}$. 
Conversely, this condition implies that $m \ge 2 + \de ( \infty ) / 2$ when $\lim_{z \to \infty } \de (z) = \de ( \infty ) \in (0, \infty )$ and $\lim_{z \to \infty } z \{ \de (z) - \de (z + 1) \} = 0$, which can be verified by considering $x^2 I_{n}^{( \de )} (x \j ^{(m)} {\j ^{(N)}}' )$ for $x \in \mathbb{N} _0$ and taking the limit as $x \to \infty $. 
(The proof is omitted.) 
In particular, in the case of the empirical Bayes estimator $\hat{\p } ^{\rm{EB}}$, the condition that $I_{n}^{( \de ^{\rm{EB}} )} ( \x ) \le 0$ for all $\x \in {\mathbb{N} _0}^{m \times N}$ implies that $m \ge 5$, while it was assumed in Corollary \ref{cor:EB} that $m \ge 7$. 
\end{remark}

\noindent
{\bf Proof of Lemma \ref{lem:B}.} \ \ Let $\x = ( x_{i, \nu } )_{1 \le i \le m, \, 1 \le \nu \le N} \in {\mathbb{N} _0}^{m \times N}$ and fix $i = 1, \dots , m$ and $\nu = 1, \dots , N$. 
The posterior mean of $1 / p_{i, \nu }$ with respect to the observation $\X = \x $ and the prior $\p \sim \prod_{{\nu }' = 1}^{N} {\rm{Dir}}_m ( \p _{{\nu }'} | a_0 , \j ^{(m)} ) \propto \prod_{{\nu }' = 1}^{N} {p_{0, {\nu }'}}^{a_0 - 1}$ is given by 
\begin{align}
E[ 1 / p_{i, \nu } | \X = \x ] &= \frac{ \int_{{D_m}^N} (1 / p_{i, \nu } ) \big\{ \prod_{{\nu }' = 1}^{N} \big( {p_{0, {\nu }'}}^{r + a_0 - 1} \prod_{j = 1}^{m} {p_{j, {\nu }'}}^{x_{j, {\nu }'}} \big) \big\} d\p }{ \int_{{D_m}^N} \big\{ \prod_{{\nu }' = 1}^{N} \big( {p_{0, {\nu }'}}^{r + a_0 - 1} \prod_{j = 1}^{m} {p_{j, {\nu }'}}^{x_{j, {\nu }'}} \big) \big\} d\p } \non \\
&= \begin{cases} \displaystyle {r + a_0 + x_{\cdot , \nu } + m - 1 \over x_{i, \nu }} & \text{if $x_{i, \nu } \ge 1$} \\ \displaystyle \infty & \text{if $x_{i, \nu } = 0$} \text{,} \end{cases} \non 
\end{align}
where $x_{\cdot , \nu } = \sum_{j = 1}^{m} x_{j, \nu }$. 
Similarly, the posterior mean of $p_{i, \nu }$ is 
\begin{align}
E[ p_{i, \nu } | \X = \x ] &= \frac{ \int_{{D_m}^N} p_{i, \nu } \big\{ \prod_{{\nu }' = 1}^{N} \big( {p_{0, {\nu }'}}^{r + a_0 - 1} \prod_{j = 1}^{m} {p_{j, {\nu }'}}^{x_{j, {\nu }'}} \big) \big\} d\p }{ \int_{{D_m}^N} \big\{ \prod_{{\nu }' = 1}^{N} \big( {p_{0, {\nu }'}}^{r + a_0 - 1} \prod_{j = 1}^{m} {p_{j, {\nu }'}}^{x_{j, {\nu }'}} \big) \big\} d\p } \non \\
&= {x_{i, \nu } + 1 \over r + a_0 + x_{\cdot , \nu } + m} < \infty \text{.} \non 
\end{align}
Therefore, for any $d \in \mathbb{R}$, the posterior expectation of the loss $(d - p_{i, \nu } )^2 / p_{i, \nu }$ can be expressed as 
\begin{align}
E[ (d - p_{i, \nu } )^2 / p_{i, \nu } | \X = \x ] &= d^2 E[ 1 / p_{i, \nu } | \X = \x ] - 2 d + E[ p_{i, \nu } | \X = \x ] \text{,} \non 
\end{align}
which is minimized at 
\begin{align}
d = {1 \over E[ 1 / p_{i, \nu } | \X = \x ]} = {x_{i, \nu } \over r + a_0 + x_{\cdot , \nu } + m - 1} \text{.} \non 
\end{align}
Hence, $\hat{\p } ^{( \de ^{( a_0 )} )} = ( X_{i, \nu } / (r + a_0 + X_{\cdot , \nu } + m - 1))_{1 \le i \le m, \, 1 \le \nu \le N}$ is a Bayes solution. 
\hfill$\Box$

\bigskip

\noindent
{\bf Proof of Lemma \ref{lem:propriety}.} \ \ Part (ii) follows immediately from part (i) since the posterior given $\X = ( \x _1 , \dots , \x _N )$ is proper for all $\x _1 , \dots , \x _N \in {\mathbb{N} _0}^m$ if and only if that given $\X = ( \bm{0} ^{(m)} , \dots , \bm{0} ^{(m)} )$, namely $\p \sim \pi _{\al , \be , g, r + a_0 , \a } ( \p )$, is proper. 
For part (i), let $J^{( \al , \be , g, a_0 , \a )} = \int_{{D_m}^N} \pi _{\al , \be , g, a_0 , \a } ( \p ) d\p $. 
Then we have 
\begin{align}
J^{( \al , \be , g, a_0 , \a )} &= \int_{0}^{\infty } t^{\al - 1} e^{- \be t} g(t) \{ B_m (t + a_0 , \a ) \} ^N dt \text{,} \non 
\end{align}
where 
\begin{align}
B_m (t + a_0 , \a ) &= \int_{D_m} \Big( {\mathring{p} _{0}}^{t + a_{0} - 1} \prod_{i = 1}^{m} {\mathring{p} _i}^{a_{i} - 1} \Big) d{\mathring{\p } _{\nu }} \non \\
&= \begin{cases} \displaystyle {\Ga (t + a_0 ) \prod_{i = 1}^{m} \Ga ( a_i ) \over \Ga (t + a_0 + a_{\cdot } )} & \text{if $t + a_0 > 0$} \\ \displaystyle \infty & \text{if $t + a_0 \le 0$} \end{cases} \non 
\end{align}
for $t \in (0, \infty )$. 
Therefore, a necessary condition for the prior to be proper is that $a_0 \ge 0$. 
Suppose that $a_0 \ge 0$. 
Then 
\begin{align}
J^{( \al , \be , g, a_0 , \a )} / \Big\{ \prod_{i = 1}^{m} \Ga ( a_i ) \Big\} ^N &= J_{1}^{( \al , \be , g, a_0 , \a )} + J_{2}^{( \al , \be , g, a_0 , \a )} \text{,} \non 
\end{align}
where 
\begin{align}
J_{1}^{( \al , \be , g, a_0 , \a )} = \int_{0}^{1} t^{\al - 1} e^{- \be t} g(t) \Big\{ {\Ga (t + a_0 ) \over \Ga (t + a_0 + a_{\cdot } )} \Big\} ^N dt \non 
\end{align}
and 
\begin{align}
J_{2}^{( \al , \be , g, a_0 , \a )} = \int_{1}^{\infty } t^{\al - 1} e^{- \be t} g(t) \Big\{ {\Ga (t + a_0 ) \over \Ga (t + a_0 + a_{\cdot } )} \Big\} ^N dt \text{.} \non 
\end{align}
The term $J_{1}^{( \al , \be , g, a_0 , \a )}$ is finite if and only if either $a_0 = 0$ and $\int_{0}^{1} t^{\al - N - 1} e^{- \be t} g(t) dt < \infty $ or $a_0 > 0$ since $\lim_{t \to 0} \Ga (t + a_0 ) / \Ga (t + a_0 + a_{\cdot } ) = \Ga ( a_0 ) / \Ga ( a_0 + a_{\cdot } )$ when $a_0 > 0$ and since $\Ga (t + 0) / \Ga (t + 0 + a_{\cdot } ) \sim t^{- 1} / \Ga ( a_{\cdot } )$ as $t \to 0$ when $a_0 = 0$. 
The term $J_{2}^{( \al , \be , g, a_0 , \a )}$ is finite if and only if $\int_{1}^{\infty } t^{\al - N a_{\cdot } - 1} e^{- \be t} g(t) dt < \infty $ since $\Ga (t + a_0 ) / \Ga (t + a_0 + a_{\cdot } ) \sim t^{- a_{\cdot }}$ as $t \to \infty $. 
This completes the proof of part (i). 
\hfill$\Box$

\bigskip

\noindent
{\bf Proof of Lemma \ref{lem:HB}.} \ \ Let $\x = ( x_{i, \nu } )_{1 \le i \le m, \, 1 \le \nu \le N} \in {\mathbb{N} _0}^{m \times N}$ and fix $i = 1, \dots , m$ and $\nu = 1, \dots , N$. 
Then it can be verified that the reciprocal of the posterior mean of $1 / p_{i, \nu }$ with respect to the observation $\X = \x $ and the prior $\p \sim \pi _{\al , \be , g, - m, \j ^{(m)}} ( \p )$ is given by 
\begin{align}
{1 \over E[ 1 / p_{i, \nu } | \X = \x ]} = \begin{cases} \displaystyle {x_{i, \nu } \over r + x_{\cdot , \nu } - 1 + \de ^{( \al , \be , g)} ( \x _{\cdot } )} & \text{if $x_{i, \nu } \ge 1$} \\ \displaystyle 0 & \text{if $x_{i, \nu } = 0$} \text{,} \end{cases} \non 
\end{align}
where $x_{\cdot , \nu } = \sum_{j = 1}^{m} x_{j, \nu }$ and $\x _{\cdot } = \big( \sum_{j = 1}^{m} x_{j, 1}, \dots , \sum_{j = 1}^{m} x_{j, N} \big) '$. 
Also, the posterior mean of $p_{i, \nu }$ is finite since $0 \le p_{i, \nu } \le 1$ and the posterior is proper. 
Therefore, for any $d \in \mathbb{R}$, the posterior expectation of the loss $(d - p_{i, \nu } )^2 / p_{i, \nu }$ can be expressed as 
\begin{align}
E[ (d - p_{i, \nu } )^2 / p_{i, \nu } | \X = \x ] &= d^2 E[ 1 / p_{i, \nu } | \X = \x ] - 2 d + E[ p_{i, \nu } | \X = \x ] \non 
\end{align}
and is minimized at $d = 1 / E[ 1 / p_{i, \nu } | \X = \x ]$. 
Hence, $\hat{\p } ^{( \al , \be , g)}$ is a Bayes solution. 
\hfill$\Box$

\bigskip

\noindent
{\bf Proof of Proposition \ref{prp:HB}.} \ \ Part (i) follows from the definition of the function $\de ^{( \al , \be , g)}$. 
Let 
\begin{align}
f_{\al , \be , g} (t) &= %
t^{\al - 1} e^{- \be t} g(t) \prod_{{\nu }' = 1}^{N} {\Ga (t + r - m) \over \Ga (t + r + z_{{\nu }'} )} \non 
\end{align}
for $t \in (0, \infty )$. 
For part (ii), suppose that $\de ^{( \al , \be , g)} ( \z ) < \infty $. 
Then by the covariance inequality we have 
\begin{align}
\de ^{( \al , \be , g)} ( \z ) / \de ^{( \al , \be , g)} ( \z + \e _{\nu }^{(N)} ) &= {\int_{0}^{\infty } t f_{\al , \be , g} (t) dt \over \int_{0}^{\infty } f_{\al , \be , g} (t) dt} / {\int_{0}^{\infty } {t \over t + r + z_{\nu }} f_{\al , \be , g} (t) dt / \int_{0}^{\infty } f_{\al , \be , g} (t) dt \over \int_{0}^{\infty } {1 \over t + r + z_{\nu }} f_{\al , \be , g} (t) dt / \int_{0}^{\infty } f_{\al , \be , g} (t) dt} \ge 1 \text{.} \non 
\end{align}
For part (iii), let $k \in \mathbb{N}$. 
Then 
\begin{align}
\de ^{( \al , \be , g)} ( \z + k \e _{\nu }^{(N)} ) &= {\int_{0}^{\infty } {t \over (t + r + z_{\nu } ) \dotsm (t + r + z_{\nu } + k - 1)} f_{\al , \be , g} (t) dt \over \int_{0}^{\infty } {1 \over (t + r + z_{\nu } ) \dotsm (t + r + z_{\nu } + k - 1)} f_{\al , \be , g} (t) dt} \text{.} \non 
\end{align}
Fix $\ep > 0$. 
Then it follows that for each $l = 0, 1$, 
\begin{align}
&\Big| \frac{ \int_{0}^{\infty } {t^l \over (t + r + z_{\nu } ) \dotsm (t + r + z_{\nu } + k - 1)} f_{\al , \be , g} (t) dt }{ \int_{0}^{\ep } {t^l \over (t + r + z_{\nu } ) \dotsm (t + r + z_{\nu } + k - 1)} f_{\al , \be , g} (t) dt } - 1 \Big| \non \\
&\le \frac{ \int_{\ep }^{\infty } {t^l \over (t + r + z_{\nu } ) \dotsm (t + r + z_{\nu } + k - 1)} f_{\al , \be , g} (t) dt }{ \int_{0}^{\ep / 2} {t^l \over (t + r + z_{\nu } ) \dotsm (t + r + z_{\nu } + k - 1)} f_{\al , \be , g} (t) dt } \non \\
&\le {(\ep / 2 + r + z_{\nu } ) \dotsm ( \ep / 2 + r + z_{\nu } + k - 1) \over ( \ep + r + z_{\nu } ) \dotsm ( \ep + r + z_{\nu } + k - 1)} \frac{ \int_{\ep }^{\infty } t^l f_{\al , \be , g} (t) dt }{ \int_{0}^{\ep / 2} t^l f_{\al , \be , g} (t) dt } \non \\
&= {\Ga ( \ep / 2 + r + z_{\nu } + k) / \Ga ( \ep / 2 + r + z_{\nu } ) \over \Ga ( \ep + r + z_{\nu } + k) / \Ga ( \ep + r + z_{\nu } )} \frac{ \int_{\ep }^{\infty } t^l f_{\al , \be , g} (t) dt }{ \int_{0}^{\ep / 2} t^l f_{\al , \be , g} (t) dt } \text{,} \non 
\end{align}
the right-hand side of which converges to zero as $k \to \infty $ since $\Ga ( \ep / 2 + r + z_{\nu } + k) / \Ga ( \ep + r + z_{\nu } + k) \sim 1 / ( \ep + r + z_{\nu } + k)^{\ep / 2}$ as $k \to \infty $. 
Therefore, 
\begin{align}
\de ^{( \al , \be , g)} ( \z + k \e _{\nu }^{(N)} ) &\sim {\int_{0}^{\ep } {t \over (t + r + z_{\nu } ) \dotsm (t + r + z_{\nu } + k - 1)} f_{\al , \be , g} (t) dt \over \int_{0}^{\ep } {1 \over (t + r + z_{\nu } ) \dotsm (t + r + z_{\nu } + k - 1)} f_{\al , \be , g} (t) dt} \le \ep \non 
\end{align}
as $k \to \infty $. 
Since $\ep $ was arbitrary, we conclude that $\lim_{\mathbb{N} \ni k \to \infty } \de ^{( \al , \be , g)} ( \z + k \e _{\nu }^{(N)} ) = 0$. 
For part (iv), let $k \in \mathbb{N} \setminus \{ 1 \} $. 
Then 
\begin{align}
&{\de ^{( \al , \be , g)} ( \z + k \j ^{(N)} ) \over 1 / \log k} \non \\
&= \frac{ \int_{0}^{\infty } ( \log k) t^{\al } e^{- \be t} g(t) \big\{ \prod_{\nu = 1}^{N} {\Ga (t + r - m) \over \Ga (t + r + z_{\nu } + k)} \big\} dt }{ \int_{0}^{\infty } t^{\al - 1} e^{- \be t} g(t) \big\{ \prod_{\nu = 1}^{N} {\Ga (t + r - m) \over \Ga (t + r + z_{\nu } + k)} \big\} dt } \non \\
&= \frac{ \int_{0}^{\infty } u^{\al } e^{- \be u / \log k} g \big( {u \over \log k} \big) \big\{ \prod_{\nu = 1}^{N} {\Ga (u / \log k + r - m) \Ga (r + z_{\nu } + k) \over \Ga (u / \log k + r + z_{\nu } + k) \Ga (r + z_{\nu } )} \big\} du }{ \int_{0}^{\infty } u^{\al - 1} e^{- \be u / \log k} g \big( {u \over \log k} \big) \big\{ \prod_{\nu = 1}^{N} {\Ga (u / \log k + r - m) \Ga (r + z_{\nu } + k) \over \Ga (u / \log k + r + z_{\nu } + k) \Ga (r + z_{\nu } )} \big\} du } \text{.} \non 
\end{align}
Now for each $l = 0, 1$ and all $u \in (0, \infty )$, we have that 
\begin{align}
&u^{\al + l - 1} e^{- \be u / \log k} g \Big( {u \over \log k} \Big) \prod_{\nu = 1}^{N} {\Ga (u / \log k + r - m) \Ga (r + z_{\nu } + k) \over \Ga (u / \log k + r + z_{\nu } + k) \Ga (r + z_{\nu } )} \non \\
&\le \frac{ \big[ \sup_{t \in (0, \infty )} \big\{ g(t) \prod_{\nu = 1}^{N} {\Ga (t + r - m) \over \Ga (t + r + z_{\nu } )} \big\} \big] u^{\al + l - 1} }{ \prod_{\nu = 1}^{N} \big\{ \big( 1 + {u / \log k \over r + z_{\nu }} \big) \dotsm \big( 1 + {u / \log k \over r + z_{\nu } + k - 1} \big) \big\} } \non \\
&\le \frac{ \big[ \sup_{t \in (0, \infty )} \big\{ g(t) \prod_{\nu = 1}^{N} {\Ga (t + r - m) \over \Ga (t + r + z_{\nu } )} \big\} \big] u^{\al + l - 1} }{ \prod_{\nu = 1}^{N} \big\{ 1 + u \big( \log {r + z_{\nu } + k \over r + z_{\nu }} \big) / \log k \big\} } \non \\
&\le \frac{ \big[ \sup_{t \in (0, \infty )} \big\{ g(t) \prod_{\nu = 1}^{N} {\Ga (t + r - m) \over \Ga (t + r + z_{\nu } )} \big\} \big] u^{\al + l - 1} }{ \prod_{\nu = 1}^{N} \big[ 1 + u \inf_{k' \in \mathbb{N} \setminus \{ 1 \} } \big\{ \big( \log {r + z_{\nu } + k' \over r + z_{\nu }} \big) / \log k' \big\} \big] } \text{,} \non 
\end{align}
where the second inequality follows since 
\begin{align}
&\Big( 1 + {u / \log k \over r + z_{\nu }} \Big) \dotsm \Big( 1 + {u / \log k \over r + z_{\nu } + k - 1} \Big) \non \\
&\ge 1 + {u \over \log k} \Big( {1 \over r + z_{\nu }} + \dots + {1 \over r + z_{\nu } + k - 1} \Big) \non \\
&\ge 1 + {u \over \log k} \log {r + z_{\nu } + k \over r + z_{\nu }} \non 
\end{align}
for every $\nu = 1, \dots , N$, and that 
\begin{align}
&\lim_{\mathbb{N} \setminus \{ 1 \} \ni k \to \infty } \Big\{ u^{\al + l - 1} e^{- \be u / \log k} g \Big( {u \over \log k} \Big) \prod_{\nu = 1}^{N} {\Ga (u / \log k + r - m) \Ga (r + z_{\nu } + k) \over \Ga (u / \log k + r + z_{\nu } + k) \Ga (r + z_{\nu } )} \Big\} \non \\
&= g(0) \Big\{ \prod_{\nu = 1}^{N} {\Ga (r - m) \over \Ga (r + z_{\nu } )} \Big\} u^{\al + l - 1} \prod_{\nu = 1}^{N} \lim_{\mathbb{N} \setminus \{ 1 \} \ni k \to \infty } {\Ga (r + z_{\nu } + k) \over \Ga (u / \log k + r + z_{\nu } + k)} \non \\
&= g(0) \Big\{ \prod_{\nu = 1}^{N} {\Ga (r - m) \over \Ga (r + z_{\nu } )} \Big\} u^{\al + l - 1} e^{- N u} \text{.} \non 
\end{align}
Thus, 
\begin{align}
\lim_{\mathbb{N} \setminus \{ 1 \} \ni k \to \infty } {\de ^{( \al , \be , g)} ( \z + k \j ^{(N)} ) \over 1 / \log k} &= \frac{ \int_{0}^{\infty } u^{\al } e^{- N u} du }{ \int_{0}^{\infty } u^{\al - 1} e^{- N u} du } = {\al \over N} \text{,} \non 
\end{align}
and the result follows. 
\hfill$\Box$

\bigskip

\noindent
{\bf Proof of Theorem \ref{thm:HB}.} \ \ First, note that $r \ge m \ge 3$ by assumption. 
Let $\De _{n}^{( \al , \be , g)} = E[ L_n ( \hat{\p } ^{( \al , \be , g)} , \p ) ] - E[ L_n ( \hat{\p } ^{\rm{U}} , \p ) ]$. 
For $\nu = 1, \dots , N$, let 
\begin{align}
\phi _{\nu }^{( \al , \be , g )} ( \X ) &= \begin{cases} \displaystyle {K( \al + 1, \be , g, r - m, \X _{\cdot } + m \j ^{(N)} ) \over K( \al , \be , g, r - m, \X _{\cdot } + m \j ^{(N)} - \e _{\nu }^{(N)} )} & \text{if $X_{\cdot , \nu } \ge 1$} \\ \displaystyle 0 & \text{if $X_{\cdot , \nu } = 0$} \end{cases} \non 
\end{align}
so that 
\begin{align}
\ph _{i, \nu }^{( \al , \be , g )} = \ph _{i, \nu }^{\rm{U}} - \ph _{i, \nu }^{\rm{U}} \phi _{\nu }^{( \al , \be , g )} ( \X ) \non 
\end{align}
for every $i = 1, \dots , m$. 
Then, by Lemma \ref{lem:hudson}, we have 
\begin{align}
\De _{n}^{( \al , \be , g)} &= E \Big[ \sum_{\nu = 1}^{n} \sum_{i = 1}^{m} \Big( {1 \over p_{i, \nu }} [ ( \ph _{i, \nu }^{\rm{U}} )^2 \{ \phi _{\nu }^{( \al , \be , g)} ( \X ) \} ^2 - 2 ( \ph _{i, \nu }^{\rm{U}} )^2 \phi _{\nu }^{( \al , \be , g)} ( \X )] + 2 \ph _{i, \nu }^{\rm{U}} \phi _{\nu }^{( \al , \be , g)} ( \X ) \Big) \Big] \non \\
&= E \Big[ \sum_{\nu = 1}^{n} \sum_{i = 1}^{m} \Big[ {X_{i, \nu } + 1 \over r + X_{\cdot , \nu }} \{ \phi _{\nu }^{( \al , \be , g)} ( \X + \e _{i, \nu }^{(m, N)} ) \} ^2 \non \\
&\quad - 2 {X_{i, \nu } + 1 \over r + X_{\cdot , \nu }} \phi _{\nu }^{( \al , \be , g)} ( \X + \e _{i, \nu }^{(m, N)} ) + 2 \ph _{i, \nu }^{\rm{U}} \phi _{\nu }^{( \al , \be , g)} ( \X ) \Big] \Big] \non \\
&= E[ I_{1, n}^{( \al , \be , g)} ( \X ) - 2 I_{2, n}^{( \al , \be , g)} ( \X ) + 2 I_{3, n}^{( \al , \be , g)} ( \X ) ] \text{,} \non 
\end{align}
where 
\begin{align}
I_{1, n}^{( \al , \be , g)} ( \x ) &= \sum_{\nu = 1}^{n} {x_{\cdot , \nu } + m \over r + x_{\cdot , \nu }} \Big\{ {K( \al + 1, \be , g, r - m, \x _{\cdot } + m \j ^{(N)} + \e _{\nu }^{(N)} ) \over K( \al , \be , g, r - m, \x _{\cdot } + m \j ^{(N)} )} \Big\} ^2 \text{,} \non \\
I_{2, n}^{( \al , \be , g)} ( \x ) &= \sum_{\nu = 1}^{n} {x_{\cdot , \nu } + m \over r + x_{\cdot , \nu }} {K( \al + 1, \be , g, r - m, \x _{\cdot } + m \j ^{(N)} + \e _{\nu }^{(N)} ) \over K( \al , \be , g, r - m, \x _{\cdot } + m \j ^{(N)} )} \text{,} \non \\
I_{3, n}^{( \al , \be , g)} ( \x ) &= \sum_{\nu = 1}^{n} {x_{\cdot , \nu } \over r + x_{\cdot , \nu } - 1} \phi _{\nu }^{( \al , \be , g )} ( \x ) \text{,} \non 
\end{align}
and $\x _{\cdot } = ( x_{\cdot , 1} , \dots , x_{\cdot , N} )' = \big( \sum_{i = 1}^{m} x_{i, 1} , \dots , \sum_{i = 1}^{m} x_{i, N} \big) '$ for $\x = %
( x_{i, \nu } )_{1 \le i \le m, \, 1 \le \nu \le N} \in {\mathbb{N} _0}^{m \times N}$. %
Since $I_{1, n}^{( \al , \be , g)} ( \bm{0} ^{(m, N)} ) - 2 I_{2, n}^{( \al , \be , g)} ( \bm{0} ^{(m, N)} ) + 2 I_{3, n}^{( \al , \be , g)} ( \bm{0} ^{(m, N)} ) < 0$, it is sufficient to show that $I_{1, n}^{( \al , \be , g)} ( \x ) - 2 I_{2, n}^{( \al , \be , g)} ( \x ) + 2 I_{3, n}^{( \al , \be , g)} ( \x ) \le 0$ for all $\x \in {\mathbb{N} _0}^{m \times N} \setminus \{ \bm{0} ^{(m, N)} \} $. 

Fix $\x = %
( x_{i, \nu } )_{1 \le i \le m, \, 1 \le \nu  \le N} \in {\mathbb{N} _0}^{m \times N} \setminus \{ \bm{0} ^{(m, N)} \} $. 
For notational simplicity, let $z _{\nu } = \sum_{i = 1}^{m} x_{i, \nu }$ for $\nu = 1, \dots , N$ and let $\z = ( z_1 , \dots , z_N )'$ and $z = \sum_{\nu = 1}^{N} z_{\nu }$. 
In addition, we use the abbreviated notation 
\begin{gather}
I_1 = I_{1, n}^{( \al , \be , g)} ( \x ) \text{,} \quad I_2 = I_{2, n}^{( \al , \be , g)} ( \x ) \text{,} \quad I_3 = I_{3, n}^{( \al , \be , g)} ( \x ) \text{,} \non \\
I = I_1 - 2 I_2 + 2 I_3 \text{,} \non \\
H(l) = {K( \al + l, \be , g, r - m, \x _{\cdot } + m \j ^{(N)} ) \over K( \al , \be , g, r - m, \x _{\cdot } + m \j ^{(N)} )} \text{,} \non \\
H(l, \pm \nu ) = {K( \al + l, \be , g, r - m, \x _{\cdot } + m \j ^{(N)} \pm \e _{\nu }^{(N)} ) \over K( \al , \be , g, r - m, \x _{\cdot } + m \j ^{(N)} )} \text{,} \non 
\end{gather}
for $l = 0, 1, 2$ and $\nu = 1, \dots , N$. 
Also, let 
\begin{align}
f_{\al , \be , g} (t) &= t^{\al - 1} e^{- \be t} g(t) \prod_{\nu = 1}^{N} {\Ga (t + r - m) \over \Ga (t + r + z_{\nu } )} \non 
\end{align}
for $t \in (0, \infty )$ so that, for example, $K( \al , \be , g, r - m, \x _{\cdot } + m \j ^{(N)} ) = \int_{0}^{\infty } f_{\al , \be , g} (t) dt$ and let 
\begin{align}
f_{\al , \be , g}^{*} (t) = {f_{\al , \be , g} (t) \over K( \al , \be , g, r - m, \x _{\cdot } + m \j ^{(N)} )} = {f_{\al , \be , g} (t) \over \int_{0}^{\infty } f_{\al , \be , g} ( t' ) d{t'}} \non 
\end{align}
for $t \in (0, \infty )$. 

For all $\nu = 1, \dots , N$ such that $z_{\nu } \neq 0$, we have that 
\begin{align}
\phi _{\nu }^{( \al , \be , g)} ( \x ) &= {H(1) \over H(0, - \nu )} = \frac{ \int_{0}^{\infty } t f_{\al , \be , g} (t) dt }{ \int_{0}^{\infty } (t + r + z_{\nu } - 1) f_{\al , \be , g} (t) dt } \non \\
&= \frac{ \int_{0}^{\infty } t f_{\al , \be , g} (t) dt }{ \int_{0}^{\infty } f_{\al , \be , g} (t) dt } \frac{ \int_{0}^{\infty } (t + r + z_{\nu } - 1) f_{\al , \be , g} (t) dt - \int_{0}^{\infty } t f_{\al , \be , g} (t) dt }{ (r + z_{\nu } - 1) \int_{0}^{\infty } (t + r + z_{\nu } - 1) f_{\al , \be , g} (t) dt } \non \\
&= {1 \over r + z_{\nu } - 1} \Big\{ H(1) - {H(1) \over H(0, - \nu )} H(1) \Big\} \non 
\end{align}
and that 
\begin{align}
{H(1) \over H(0, - \nu )} / H(1, \nu ) &= \frac{ \int_{0}^{\infty } t f_{\al , \be , g}^{*} (t) dt }{ \int_{0}^{\infty } (t + r + z_{\nu } - 1) f_{\al , \be , g}^{*} (t) dt \int_{0}^{\infty } {t \over t + r + z_{\nu }} f_{\al , \be , g}^{*} (t) dt } \non \\
&\ge \frac{ \int_{0}^{\infty } t f_{\al , \be , g}^{*} (t) dt }{ \int_{0}^{\infty } t {t + r + z_{\nu } - 1 \over t + r + z_{\nu }} f_{\al , \be , g}^{*} (t) dt } \ge 1 \non 
\end{align}
by the covariance inequality. 
Therefore, 
\begin{align}
I_3 &\le \sum_{\nu = 1}^{n} {z_{\nu } \over (r + z_{\nu } - 1)^2} \{ H(1) - H(1, \nu ) H(1) \} \non \\
&\le \sum_{\nu = 1}^{n} {z_{\nu } + 2 \over (r + z_{\nu } )^2} \{ H(1) - H(1, \nu ) H(1) \} \non \\
&= \sum_{\nu = 1}^{n} {z_{\nu } + 2 \over (r + z_{\nu } )^2} H(1) - \sum_{\nu = 1}^{n} {z_{\nu } + 2 \over (r + z_{\nu } )^2} H(1, \nu ) H(1) \text{.} \label{tHBp1} 
\end{align}
Since 
\begin{align}
H(1, \nu ) &= \int_{0}^{\infty } {t \over t + r + z_{\nu }} f_{\al , \be , g}^{*} (t) dt \non \\
&= \int_{0}^{\infty } {t \over r + z_{\nu }} \Big( 1 - {t \over t + r + z_{\nu }} \Big) f_{\al , \be , g}^{*} (t) dt \non \\
&= {1 \over r + z_{\nu }} \{ H(1) - H(2, \nu ) \} \text{,} \non 
\end{align}
for all $\nu = 1, \dots , N$, it follows that 
\begin{align}
I_2 &= \sum_{\nu = 1}^{n} {z_{\nu } + m \over r + z_{ \nu }} H(1, \nu ) \non \\
&= \sum_{\nu = 1}^{n} {z_{\nu } + m \over (r + z_{ \nu } )^2} H(1) - \sum_{\nu = 1}^{n} {z_{\nu } + m \over (r + z_{ \nu } )^2} H(2, \nu ) \text{.} \label{tHBp2} 
\end{align}
Now, by the covariance inequality, 
\begin{align}
\sum_{\nu = 1}^{n} {z_{\nu } + m \over (r + z_{ \nu } )^2} H(2, \nu ) &\le {1 \over n} \Big\{ \sum_{\nu = 1}^{n} {1 \over (r + z_{ \nu } )^2} \Big\} \sum_{\nu = 1}^{n} ( z_{\nu } + m) H(2, \nu ) \text{.} \label{tHBp3} 
\end{align}
By integration by parts, 
\begin{align}
\infty &> ( \al + 1) \int_{0}^{\infty } t f_{\al , \be , g} (t) dt \non \\
&= \int_{0}^{\infty } ( \al + 1) t^{\al } e^{- \be t} g(t) \Big\{ \prod_{{\nu }' = 1}^{N} {\Ga (t + r - m) \over \Ga (t + r + z_{{\nu }'} )} \Big\} dt \non \\
&= \lim_{\ep \to 0} \Big( \Big[ t^{\al + 1} e^{- \be t} g(t) \Big\{ \prod_{{\nu }' = 1}^{N} {\Ga (t + r - m) \over \Ga (t + r + z_{{\nu }'} )} \Big\} \Big] _{\ep }^{1 / \ep } \non \\
&\quad - \int_{\ep }^{1 / \ep } t^{\al + 1} \Big[ {\partial \over \partial t} \Big\{ e^{- \be t} g(t) \prod_{{\nu }' = 1}^{N} {\Ga (t + r - m) \over \Ga (t + r + z_{{\nu }'} )} \Big\} \Big] dt \Big) \non \\
&= \int_{0}^{\infty } t^{\al + 1} e^{- \be t} g(t) \Big\{ \prod_{{\nu }' = 1}^{N} {\Ga (t + r - m) \over \Ga (t + r + z_{{\nu }'} )} \Big\} \Big\{ \be + {- g' (t) \over g(t)} \Big\} dt \non \\
&\quad + \sum_{\nu = 1}^{N} \sum_{k = 1}^{z_{\nu } + m} \int_{0}^{\infty } t^{\al + 1} e^{- \be t} g(t) \Big\{ \prod_{{\nu }' = 1}^{N} {\Ga (t + r - m) \over \Ga (t + r + z_{{\nu }'} )} \Big\} {1 \over t + r - m + k - 1} dt \text{,} \non 
\end{align}
where the last equality follows from the assumptions of the theorem since $\Ga (t) \sim t^{- 1}$ as $t \to 0$ while $\prod_{{\nu }' = 1}^{N} \{ \Ga (t + r - m) / \Ga (t + r + z_{{\nu }'} ) \} \sim t^{- z - N m}$ as $t \to \infty $ and since $- g' (t) \ge 0$ and $|t / (t + r - m + k - 1)| \le 1$ for all $t \in (0, \infty )$ and $k = 1, 2, \dotsc $. 
Therefore, 
\begin{align}
&\sum_{\nu  = 1}^{n} ( z_{\nu } + m) H(2, \nu ) \non \\
&= \sum_{\nu = 1}^{n} \sum_{k = 1}^{z_{\nu } + m} \int_{0}^{\infty } t^{\al + 1} e^{- \be t} g(t) \Big\{ \prod_{{\nu }' = 1}^{N} {\Ga (t + r - m) \over \Ga (t + r + z_{{\nu }'} )} \Big\} {1 \over t + r + z_{\nu }} dt / \int_{0}^{\infty } f_{\al , \be , g} (t) dt \non \\
&\le \sum_{\nu = 1}^{N} \sum_{k = 1}^{z_{\nu } + m} \int_{0}^{\infty } t^{\al + 1} e^{- \be t} g(t) \Big\{ \prod_{{\nu }' = 1}^{N} {\Ga (t + r - m) \over \Ga (t + r + z_{{\nu }'} )} \Big\} {1 \over t + r - m + k - 1} dt / \int_{0}^{\infty } f_{\al , \be , g} (t) dt \non \\
&\le \Big[ ( \al + 1) \int_{0}^{\infty } t f_{\al , \be , g} (t) dt - \be \int_{0}^{\infty } t^{\al + 1} e^{- \be t} g(t) \Big\{ \prod_{{\nu }' = 1}^{N} {\Ga (t + r - m) \over \Ga (t + r + z_{{\nu }'} )} \Big\} dt \Big] / \int_{0}^{\infty } f_{\al , \be , g} (t) dt \non \\
&= ( \al + 1) H(1) - \be H(2) \le ( \al + 1) H(1) - \be H(1) H(1) \text{,} \label{tHBp4} 
\end{align}
where the last inequality follows since 
\begin{align}
H(2) &= \int_{0}^{\infty } t^2 f_{\al , \be , g}^{*} (t) dt \ge \Big\{ \int_{0}^{\infty } t f_{\al , \be , g}^{*} (t) dt \Big\} ^2 = \{ H(1) \} ^2 \non 
\end{align}
by the covariance inequality. 
Combining (\ref{tHBp2}), (\ref{tHBp3}), and (\ref{tHBp4}) gives 
\begin{align}
I_2 &\ge \sum_{\nu = 1}^{n} {z_{\nu } + m \over (r + z_{ \nu } )^2} H(1) - {1 \over n} \Big\{ \sum_{\nu = 1}^{n} {1 \over (r + z_{ \nu } )^2} \Big\} \sum_{\nu = 1}^{n} ( z_{\nu } + m) H(2, \nu ) \non \\
&\ge \sum_{\nu = 1}^{n} {z_{\nu } + m - ( \al + 1) / n \over (r + z_{ \nu } )^2} H(1) + {\be \over n} \sum_{\nu = 1}^{n} {1 \over (r + z_{ \nu } )^2} H(1) H(1) \non \\
&\ge \sum_{\nu = 1}^{n} {z_{\nu } + m - ( \al + 1) / n \over (r + z_{ \nu } )^2} H(1) + {\be \over n} \sum_{\nu = 1}^{n} {r + z_{\nu } \over (r + z_{ \nu } )^2} H(1, \nu ) H(1) \label{tHBp5} 
\end{align}
since 
\begin{align}
H(1) &= \int_{0}^{\infty } t f_{\al , \be , g}^{*} (t) dt \ge \int_{0}^{\infty } {r + z_{\nu } \over t + r + z_{\nu }} t f_{\al , \be , g}^{*} (t) dt = (r + z_{\nu } ) H(1, \nu ) \non 
\end{align}
for all $\nu = 1, \dots , N$. 
Finally, 
\begin{align}
I_1 &= \sum_{\nu = 1}^{n} {z_{\nu } + m \over r + z_{\nu }} \{ H(1, \nu ) \} ^2 \non \\
&\le \sum_{\nu = 1}^{n} {z_{\nu } + m \over (r + z_{\nu } )^2} H(1, \nu ) H(1) \text{.} \label{tHBp6} 
\end{align}
Also 
\begin{align}
H(1, \nu ) &= \int_{0}^{\infty } {t \over t + r + z_{\nu }} f_{\al , \be , g}^{*} (t) dt \le \int_{0}^{\infty } f_{\al , \be , g}^{*} (t) dt = 1 \non 
\end{align}
for all $\nu = 1, \dots , N$. 
Hence, combining (\ref{tHBp1}), (\ref{tHBp5}), and (\ref{tHBp6}), we obtain 
\begin{align}
I %
&\le \sum_{\nu = 1}^{n} {- z_{\nu } + m - 4 \over (r + z_{\nu } )^2} H(1, \nu ) H(1) \non \\
&\quad - 2 \sum_{\nu = 1}^{n} {m - 2 - ( \al + 1) / n \over (r + z_{ \nu } )^2} H(1) - 2 {\be \over n} \sum_{\nu = 1}^{n} {r + z_{\nu } \over (r + z_{ \nu } )^2} H(1, \nu ) H(1) \non \\
&\le \sum_{\nu = 1}^{n} {- z_{\nu } - m + 2 ( \al + 1) / n - 2 ( \be / n) (r + z_{\nu } ) \over (r + z_{\nu } )^2} H(1, \nu ) H(1) \non \\
&\le 0 \text{,} \non 
\end{align}
where the second and third inequalities follow from (\ref{eq:assumption_HB}), and this completes the proof. 
\hfill$\Box$

\bigskip

\noindent
{\bf Proof of Lemma \ref{lem:Jeff}.} \ \ Let $\mathring{\X } = ( \mathring{X} _1 , \dots , \mathring{X} _m )' \sim {\rm{NM}}_m (r, \mathring{\p } )$. 
Then the square root of the determinant of the information matrix corresponding to %
this distribution is 
\begin{align}
&\sqrt{\Big| \Big( E \Big[ - {\partial ^2 \over \partial \mathring{p} _i \partial \mathring{p} _j} \log {\rm{NM}}_m ( \mathring{\X } | r, \mathring{p} ) \Big] \Big) _{1 \le i, j \le m} \Big| } = \sqrt{\Big| \Big( E \Big[ {r \over {{\mathring{p}}_0}^2 } + \de _{i, j}^{(m)} {\mathring{X} _i \over {{\mathring{p}} _i}^2} \Big] \Big) _{1 \le i, j \le m} \Big| } \non \\
&= \sqrt{| \D ( \mathring{\p } ) + (r / {\mathring{p} _0}^2 ) \j ^{(m)} {\j ^{(m)}}' |} = \sqrt{| \D ( \mathring{\p } )| [1 + (r / {\mathring{p} _0}^2 ) {\j ^{(m)}}' \{ \D ( \mathring{\p } ) \} ^{- 1} \j ^{(m)} ]} \\
&= \sqrt{{r^m \over {\mathring{p} _0}^m} \Big( \prod_{i = 1}^{m} {1 \over \mathring{p} _i} \Big) \Big( 1 + {{\mathring{p}} _{\cdot } \over \mathring{p} _0} \Big) } \propto {\rm{Dir}}_m \Big( \mathring{\p } \Big| {1 - m \over 2}, {1 \over 2} \j ^{(m)} \Big) \text{,} \non 
\end{align}
where $\D ( \p ) = (r / \mathring{p} _0 ) {\boldmath{\diag }} (1 / \mathring{p} _1 , \dots , 1 / \mathring{p} _m )$. 
This is the desired result. 
\hfill$\Box$

\bigskip

\noindent
{\bf Proof of Lemma \ref{lem:KL}.} \ \ Let $\x = ( x_{i, \nu } )_{1 \le i \le m, \, 1 \le \nu \le N} \in {\mathbb{N} _0}^{m \times N}$ and fix $i = 1, \dots , m$ and $\nu = 1, \dots , N$. Since the prior density is strictly positive and the posterior is proper, the posterior mean of $p_{i, \nu }$ with respect to the observation $\X = \x $, denoted $E[ p_{i, \nu } | \X = \x ]$, satisfies $E[ p_{i, \nu } | \X = \x ] \in (0, \infty )$. 
Also, $E[ p_{i, \nu } \log (1 / p_{i, \nu } ) | \X = \x ] \in (0, \infty )$. 
Therefore, for any $\tilde{d} \in (0, \infty )$, the posterior expectation of the loss $\tilde{d} - p_{i, \nu } - p_{i, \nu } \log ( \tilde{d} / p_{i, \nu } )$ can be expressed as 
\begin{align}
&E[ \tilde{d} - p_{i, \nu } - p_{i, \nu } \log ( \tilde{d} / p_{i, \nu } ) | \X = \x ] \non \\
&= \tilde{d} - E[ p_{i, \nu } | \X = \x ] \log \tilde{d} - E[ p_{i, \nu } | \X = \x ] - E[ p_{i, \nu } \log (1 / p_{i, \nu } ) | \X = \x ] \non 
\end{align}
and thus is minimized at $\tilde{d} = E[ p_{i, \nu } | \X = \x ]$, which yields the desired result. 
\hfill$\Box$

\bigskip

\noindent
{\bf Proof of Proposition \ref{prp:KL}.} \ \ The proof is similar to that of Proposition \ref{prp:HB}. 
Part (i) follows from the definition. 
Let 
\begin{align}
f_{\al , \be , g, a_0 , \a } (t) &= %
t^{\al - 1} e^{- \be t} g(t) \prod_{{\nu }'' = 1}^{N} {\Ga (t + r + a_0 ) \over \Ga (t + r + a_0 + z_{{\nu }''} + a_{\cdot } )} \non 
\end{align}
for $t \in (0, \infty )$ so that $K( \al , \be , g, r + a_0 , \z + a_{\cdot } \j ^{(N)} ) = \int_{0}^{\infty } f_{\al , \be , g, a_0 , \a } (t) dt$. 
Then part (ii) follows since 
\begin{align}
\frac{ \de _{\nu }^{( \al , \be , g, a_0 , \a )} ( \z ) }{ \de _{\nu }^{( \al , \be , g, a_0 , \a )} ( \z + \e _{{\nu }'}^{(N)} ) } &= \frac{ \displaystyle {\int_{0}^{\infty } t {f_{\al , \be , g, a_0 , \a } (t) \over t + r + a_0 + z_{\nu } + a_{\cdot }} dt \over \int_{0}^{\infty } {f_{\al , \be , g, a_0 , \a } (t) \over t + r + a_0 + z_{\nu } + a_{\cdot }} dt} }{ \displaystyle {\int_{0}^{\infty } {t \over t + r + a_0 + z_{{\nu }'} + a_{\cdot } + \de _{\nu , {\nu }'}^{(N)}} {f_{\al , \be , g, a_0 , \a } (t) \over t + r + a_0 + z_{\nu } + a_{\cdot }} dt / \int_{0}^{\infty } {f_{\al , \be , g, a_0 , \a } (t) \over t + r + a_0 + z_{\nu } + a_{\cdot }} dt \over \int_{0}^{\infty } {1 \over t + r + a_0 + z_{{\nu }'} + a_{\cdot } + \de _{\nu , {\nu }'}^{(N)}} {f_{\al , \be , g, a_0 , \a } (t) \over t + r + a_0 + z_{\nu } + a_{\cdot }} dt / \int_{0}^{\infty } {f_{\al , \be , g, a_0 , \a } (t) \over t + r + a_0 + z_{\nu } + a_{\cdot }} dt} } \ge 1 \non 
\end{align}
by the covariance inequality. 
For part (iii), let $k \in \mathbb{N}$. 
Then 
\begin{align}
\de _{\nu }^{( \al , \be , g, a_0 , \a )} ( \z + k \e _{{\nu }'}^{(N)} ) &= {\int_{0}^{\infty } {t \over (t + r + a_0 + z_{{\nu }'} + a_{\cdot } + \de _{\nu , {\nu }'}^{(N)} ) \dotsm (t + r + a_0 + z_{{\nu }'} + a_{\cdot } + \de _{\nu , {\nu }'}^{(N)} + k - 1)} {f_{\al , \be , g, a_0 , \a } (t) \over t + r + a_0 + z_{\nu } + a_{\cdot }} dt \over \int_{0}^{\infty } {1 \over (t + r + a_0 + z_{{\nu }'} + a_{\cdot } + \de _{\nu , {\nu }'}^{(N)} ) \dotsm (t + r + a_0 + z_{{\nu }'} + a_{\cdot } + \de _{\nu , {\nu }'}^{(N)} + k - 1)} {f_{\al , \be , g, a_0 , \a } (t) \over t + r + a_0 + z_{\nu } + a_{\cdot }} dt} \text{.} \non 
\end{align}
Fix $\ep > 0$. 
Then it follows that for each $l = 0, 1$, 
\begin{align}
&\Big| \frac{ \int_{0}^{\infty } {t^l \over (t + r + a_0 + z_{{\nu }'} + a_{\cdot } + \de _{\nu , {\nu }'}^{(N)} ) \dotsm (t + r + a_0 + z_{{\nu }'} + a_{\cdot } + \de _{\nu , {\nu }'}^{(N)} + k - 1)} {f_{\al , \be , g, a_0 , \a } (t) \over t + r + a_0 + z_{\nu } + a_{\cdot }} dt }{ \int_{0}^{\ep } {t^l \over (t + r + a_0 + z_{{\nu }'} + a_{\cdot } + \de _{\nu , {\nu }'}^{(N)} ) \dotsm (t + r + a_0 + z_{{\nu }'} + a_{\cdot } + \de _{\nu , {\nu }'}^{(N)} + k - 1)} {f_{\al , \be , g, a_0 , \a } (t) \over t + r + a_0 + z_{\nu } + a_{\cdot }} dt } - 1 \Big| \non \\
&\le \frac{ \int_{\ep }^{\infty } {t^l \over (t + r + a_0 + z_{{\nu }'} + a_{\cdot } + \de _{\nu , {\nu }'}^{(N)} ) \dotsm (t + r + a_0 + z_{{\nu }'} + a_{\cdot } + \de _{\nu , {\nu }'}^{(N)} + k - 1)} {f_{\al , \be , g, a_0 , \a } (t) \over t + r + a_0 + z_{\nu } + a_{\cdot }} dt }{ \int_{0}^{\ep / 2} {t^l \over (t + r + a_0 + z_{{\nu }'} + a_{\cdot } + \de _{\nu , {\nu }'}^{(N)} ) \dotsm (t + r + a_0 + z_{{\nu }'} + a_{\cdot } + \de _{\nu , {\nu }'}^{(N)} + k - 1)} {f_{\al , \be , g, a_0 , \a } (t) \over t + r + a_0 + z_{\nu } + a_{\cdot }} dt } \non \\
&\le {( \ep / 2 + r + a_0 + z_{{\nu }'} + a_{\cdot } + \de _{\nu , {\nu }'}^{(N)} ) \dotsm ( \ep / 2 + r + a_0 + z_{{\nu }'} + a_{\cdot } + \de _{\nu , {\nu }'}^{(N)} + k - 1) \over ( \ep + r + a_0 + z_{{\nu }'} + a_{\cdot } + \de _{\nu , {\nu }'}^{(N)} ) \dotsm ( \ep + r + a_0 + z_{{\nu }'} + a_{\cdot } + \de _{\nu , {\nu }'}^{(N)} + k - 1)} \frac{ \int_{\ep }^{\infty } t^l {f_{\al , \be , g, a_0 , \a } (t) \over t + r + a_0 + z_{\nu } + a_{\cdot }} dt }{ \int_{0}^{\ep / 2} t^l {f_{\al , \be , g, a_0 , \a } (t) \over t + r + a_0 + z_{\nu } + a_{\cdot }} dt } \non \\
&= {\Ga ( \ep / 2 + r + a_0 + z_{{\nu }'} + a_{\cdot } + \de _{\nu , {\nu }'}^{(N)} + k) / \Ga ( \ep / 2 + r + a_0 + z_{{\nu }'} + a_{\cdot } + \de _{\nu , {\nu }'}^{(N)} ) \over \Ga ( \ep + r + a_0 + z_{{\nu }'} + a_{\cdot } + \de _{\nu , {\nu }'}^{(N)} + k) / \Ga ( \ep + r + a_0 + z_{{\nu }'} + a_{\cdot } + \de _{\nu , {\nu }'}^{(N)} )} \frac{ \int_{\ep }^{\infty } t^l {f_{\al , \be , g, a_0 , \a } (t) \over t + r + a_0 + z_{\nu } + a_{\cdot }} dt }{ \int_{0}^{\ep / 2} t^l {f_{\al , \be , g, a_0 , \a } (t) \over t + r + a_0 + z_{\nu } + a_{\cdot }} dt } \text{,} \non 
\end{align}
the right-hand side of which converges to zero as $k \to \infty $ since $\Ga ( \ep / 2 + r + a_0 + z_{{\nu }'} + a_{\cdot } + \de _{\nu , {\nu }'}^{(N)} + k) / \Ga ( \ep + r + a_0 + z_{{\nu }'} + a_{\cdot } + \de _{\nu , {\nu }'}^{(N)} + k) \sim 1 / ( \ep + r + a_0 + z_{{\nu }'} + a_{\cdot } + \de _{\nu , {\nu }'}^{(N)} + k)^{\ep / 2}$ as $k \to \infty $. 
Therefore, 
\begin{align}
\de _{\nu }^{( \al , \be , g, a_0 , \a )} ( \z + k \e _{{\nu }'}^{(N)} ) &\sim {\int_{0}^{\ep } {t \over (t + r + a_0 + z_{{\nu }'} + a_{\cdot } + \de _{\nu , {\nu }'}^{(N)} ) \dotsm (t + r + a_0 + z_{{\nu }'} + a_{\cdot } + \de _{\nu , {\nu }'}^{(N)} + k - 1)} {f_{\al , \be , g, a_0 , \a } (t) \over t + r + a_0 + z_{\nu } + a_{\cdot }} dt \over \int_{0}^{\ep } {1 \over (t + r + a_0 + z_{{\nu }'} + a_{\cdot } + \de _{\nu , {\nu }'}^{(N)} ) \dotsm (t + r + a_0 + z_{{\nu }'} + a_{\cdot } + \de _{\nu , {\nu }'}^{(N)} + k - 1)} {f_{\al , \be , g, a_0 , \a } (t) \over t + r + a_0 + z_{\nu } + a_{\cdot }} dt} \le \ep \non 
\end{align}
as $k \to \infty $. 
Since $\ep $ was arbitrary, we conclude that $\lim_{\mathbb{N} \ni k \to \infty } \de _{\nu }^{( \al , \be , g, a_0 , \a )} ( \z + k \e _{{\nu }'}^{(N)} ) = 0$. 
For part (iv), let $k \in \mathbb{N} \setminus \{ 1 \} $. 
Then 
\begin{align}
&{\de _{\nu }^{( \al , \be , g, a_0 , \a )} ( \z + k \j ^{(N)} ) \over 1 / \log k} \non \\
&= \frac{ \int_{0}^{\infty } ( \log k) t^{\al } e^{- \be t} g(t) \big\{ \prod_{{\nu }' = 1}^{N} {\Ga (t + r + a_0 ) \over \Ga (t + r + a_0 + z_{{\nu }'} + a_{\cdot } + \de _{\nu , {\nu }'}^{(N)} + k)} \big\} dt }{ \int_{0}^{\infty } t^{\al - 1} e^{- \be t} g(t) \big\{ \prod_{{\nu }' = 1}^{N} {\Ga (t + r + a_0 ) \over \Ga (t + r + a_0 + z_{{\nu }'} + a_{\cdot } + \de _{\nu , {\nu }'}^{(N)} + k)} \big\} dt } \non \\
&= \frac{ \int_{0}^{\infty } u^{\al } e^{- \be u / \log k} g \big( {u \over \log k} \big) \big\{ \prod_{{\nu }' = 1}^{N} {\Ga (u / \log k + r + a_0 ) \Ga (r + a_0 + z_{{\nu }'} + a_{\cdot } + \de _{\nu , {\nu }'}^{(N)} + k) \over \Ga (u / \log k + r + a_0 + z_{{\nu }'} + a_{\cdot } + \de _{\nu , {\nu }'}^{(N)} + k) \Ga (r + a_0 + z_{{\nu }'} + a_{\cdot } + \de _{\nu , {\nu }'}^{(N)} )} \big\} du }{ \int_{0}^{\infty } u^{\al - 1} e^{- \be u / \log k} g \big( {u \over \log k} \big) \big\{ \prod_{{\nu }' = 1}^{N} {\Ga (u / \log k + r + a_0 ) \Ga (r + a_0 + z_{{\nu }'} + a_{\cdot } + \de _{\nu , {\nu }'}^{(N)} + k) \over \Ga (u / \log k + r + a_0 + z_{{\nu }'} + a_{\cdot } + \de _{\nu , {\nu }'}^{(N)} + k) \Ga (r + a_0 + z_{{\nu }'} + a_{\cdot } + \de _{\nu , {\nu }'}^{(N)} )} \big\} du } \text{.} \non 
\end{align}
Now for each $l = 0, 1$ and all $u \in (0, \infty )$, we have that 
\begin{align}
&u^{\al + l - 1} e^{- \be u / \log k} g \Big( {u \over \log k} \Big) \non \\
&\quad \times \prod_{{\nu }' = 1}^{N} {\Ga (u / \log k + r + a_0 ) \Ga (r + a_0 + z_{{\nu }'} + a_{\cdot } + \de _{\nu , {\nu }'}^{(N)} + k) \over \Ga (u / \log k + r + a_0 + z_{{\nu }'} + a_{\cdot } + \de _{\nu , {\nu }'}^{(N)} + k) \Ga (r + a_0 + z_{{\nu }'} + a_{\cdot } + \de _{\nu , {\nu }'}^{(N)} )} \non \\
&\le \frac{ \big[ \sup_{t \in (0, \infty )} \big\{ g(t) \prod_{{\nu }' = 1}^{N} {\Ga (t + r + a_0 ) \over \Ga (t + r + a_0 + z_{{\nu }'} + a_{\cdot } + \de _{\nu , {\nu }'}^{(N)} )} \big\} \big] u^{\al + l - 1} }{ \prod_{{\nu }' = 1}^{N} \big\{ \big( 1 + {u / \log k \over r + a_0 + z_{{\nu }'} + a_{\cdot } + \de _{\nu , {\nu }'}^{(N)}} \big) \dotsm \big( 1 + {u / \log k \over r + a_0 + z_{{\nu }'} + a_{\cdot } + \de _{\nu , {\nu }'}^{(N)} + k - 1} \big) \big\} } \non \\
&\le \frac{ \big[ \sup_{t \in (0, \infty )} \big\{ g(t) \prod_{{\nu }' = 1}^{N} {\Ga (t + r + a_0 ) \over \Ga (t + r + a_0 + z_{{\nu }'} + a_{\cdot } + \de _{\nu , {\nu }'}^{(N)} )} \big\} \big] u^{\al + l - 1} }{ \prod_{{\nu }' = 1}^{N} \big\{ 1 + u \big( \log {r + a_0 + z_{{\nu }'} + a_{\cdot } + \de _{\nu , {\nu }'}^{(N)} + k \over r + a_0 + z_{{\nu }'} + a_{\cdot } + \de _{\nu , {\nu }'}^{(N)}} \big) / \log k \big\} } \non \\
&\le \frac{ \big[ \sup_{t \in (0, \infty )} \big\{ g(t) \prod_{{\nu }' = 1}^{N} {\Ga (t + r + a_0 ) \over \Ga (t + r + a_0 + z_{{\nu }'} + a_{\cdot } + \de _{\nu , {\nu }'}^{(N)} )} \big\} \big] u^{\al + l - 1} }{ \prod_{{\nu }' = 1}^{N} \big[ 1 + u \inf _{k' \in \mathbb{N} \setminus \{ 1 \} } \big\{ \big( \log {r + a_0 + z_{{\nu }'} + a_{\cdot } + \de _{\nu , {\nu }'}^{(N)} + k' \over r + a_0 + z_{{\nu }'} + a_{\cdot } + \de _{\nu , {\nu }'}^{(N)}} \big) / \log k' \big\} \big] } \text{,} \non 
\end{align}
where the second inequality follows since 
\begin{align}
&\Big( 1 + {u / \log k \over r + a_0 + z_{{\nu }'} + a_{\cdot } + \de _{\nu , {\nu }'}^{(N)}} \Big) \dotsm \Big( 1 + {u / \log k \over r + a_0 + z_{{\nu }'} + a_{\cdot } + \de _{\nu , {\nu }'}^{(N)} + k - 1} \Big) \non \\
&\ge 1 + {u \over \log k} \Big( {1 \over r + a_0 + z_{{\nu }'} + a_{\cdot } + \de _{\nu , {\nu }'}^{(N)}} + \dots + {1 \over r + a_0 + z_{{\nu }'} + a_{\cdot } + \de _{\nu , {\nu }'}^{(N)} + k - 1} \Big) \non \\
&\ge 1 + {u \over \log k} \log {r + a_0 + z_{{\nu }'} + a_{\cdot } + \de _{\nu , {\nu }'}^{(N)} + k \over r + a_0 + z_{{\nu }'} + a_{\cdot } + \de _{\nu , {\nu }'}^{(N)}} \non 
\end{align}
for every ${\nu }' = 1, \dots , N$, and that 
\begin{align}
&\lim_{\mathbb{N} \setminus \{ 1 \} \ni k \to \infty } \Big\{ u^{\al + l - 1} e^{- \be u / \log k} g \Big( {u \over \log k} \Big) \non \\
&\quad \times \prod_{{\nu }' = 1}^{N} {\Ga (u / \log k + r + a_0 ) \Ga (r + a_0 + z_{{\nu }'} + a_{\cdot } + \de _{\nu , {\nu }'}^{(N)} + k) \over \Ga (u / \log k + r + a_0 + z_{{\nu }'} + a_{\cdot } + \de _{\nu , {\nu }'}^{(N)} + k) \Ga (r + a_0 + z_{{\nu }'} + a_{\cdot } + \de _{\nu , {\nu }'}^{(N)} )} \Big\} \non \\
&= g(0) \Big\{ \prod_{{\nu }' = 1}^{N} {\Ga (r + a_0 ) \over \Ga (r + a_0 + z_{{\nu }'} + a_{\cdot } + \de _{\nu , {\nu }'}^{(N)} )} \Big\} u^{\al + l - 1} \non \\
&\quad \times \prod_{{\nu }' = 1}^{N} \lim_{\mathbb{N} \setminus \{ 1 \} \ni k \to \infty } {\Ga (r + a_0 + z_{{\nu }'} + a_{\cdot } + \de _{\nu , {\nu }'}^{(N)} + k) \over \Ga (u / \log k + r + a_0 + z_{{\nu }'} + a_{\cdot } + \de _{\nu , {\nu }'}^{(N)} + k)} \non \\
&= g(0) \Big\{ \prod_{{\nu }' = 1}^{N} {\Ga (r + a_0 ) \over \Ga (r + a_0 + z_{{\nu }'} + a_{\cdot } + \de _{\nu , {\nu }'}^{(N)} )} \Big\} u^{\al + l - 1} e^{- N u} \text{.} \non 
\end{align}
Thus, 
\begin{align}
\lim_{\mathbb{N} \setminus \{ 1 \} \ni k \to \infty } {\de _{\nu }^{( \al , \be , g, a_0 , \a )} ( \z + k \j ^{(N)} ) \over 1 / \log k} &= \frac{ \int_{0}^{\infty } u^{\al } e^{- N u} du }{ \int_{0}^{\infty } u^{\al - 1} e^{- N u} du } = {\al \over N} \text{,} \non 
\end{align}
and the result follows. 
\hfill$\Box$

\bigskip

The following lemma will be used in the proof of Theorem \ref{thm:KL}. 

\begin{lem}
\label{lem:digamma} 
Let $u_1 , u_2 > 0$. 
Then 
\begin{align}
{{\Ga }' ( u_1 + u_2 ) \over \Ga ( u_1 + u_2 )} - {{\Ga }' ( u_1 ) \over \Ga ( u_1 )} \ge {u_2 \over u_1 + u_2} \text{.} \non 
\end{align}
\end{lem}

\noindent
{\bf Proof%
.} \ \ We have 
\begin{align}
1 &= u_1 \sum_{k = 0}^{\infty } \Big( {1 \over u_1 + k} - {1 \over u_1 + k + 1} \Big) = \sum_{k = 0}^{\infty } {u_1 \over ( u_1 + k) ( u_1 + k + 1)} \non \\
&\le \sum_{k = 0}^{\infty } {u_1 \over ( u_1 + k)^2} \le \sum_{k = 0}^{\infty } {u_1 + u_2 \over ( u_1 + k) ( u_1 + u_2 + k)} \non \\
&= {u_1 + u_2 \over u_2} \Big\{ {{\Ga }' ( u_1 + u_2 ) \over \Ga ( u_1 + u_2 )} - {{\Ga }' ( u_1 ) \over \Ga ( u_1 )} \Big\} \text{,} \non 
\end{align}
which is the desired result. 
\hfill$\Box$

\bigskip

\noindent
{\bf Proof of Theorem \ref{thm:KL}.} \ \ The proof is similar to that of Theorem \ref{thm:HB}. 
First, note that $r > 2$ and $a_{\cdot } > 1$ by assumption. 
Let $\De _{n}^{( \al , \be , g, a_0 , \a )} = E[ \tilde{L} _n ( \hat{\p } ^{( \al , \be , g, a_0 , \a )} , \p ) ] - E[ \tilde{L} _n ( \hat{\p } ^{( a_0 , \a )} , \p ) ]$. 
For $\nu = 1, \dots , N$, let 
\begin{align}
\phi _{\nu }^{( \al , \be , g, a_0 , \a )} ( \X ) &= {K( \al + 1, \be , g, r + a_0 , \X _{\cdot } + a_{\cdot } \j ^{(N)} + \e _{\nu }^{(N)} ) \over K( \al , \be , g, r + a_0 , \X _{\cdot } + a_{\cdot } \j ^{(N)} )} \non 
\end{align}
so that 
\begin{align}
\ph _{i, \nu }^{( \al , \be , g, a_0 , \a )} = \ph _{i, \nu }^{( a_0 , \a )} - \ph _{i, \nu }^{( a_0 , \a )} \phi _{\nu }^{( \al , \be , g, a_0 , \a )} ( \X ) \non 
\end{align}
for every $i = 1, \dots , m$. 
By Lemma \ref{lem:hudson}, we have 
\begin{align}
\De _{n}^{( \al , \be , g, a_0 , \a )} &= E \Big[ \sum_{\nu = 1}^{n} \sum_{i = 1}^{m} \Big[ - \ph _{i, \nu }^{( a_0 , \a )} \phi _{\nu }^{( \al , \be , g, a_0 , \a )} ( \X ) + p_{i, \nu } \log \Big\{ 1 + {\de _{\nu }^{( \al , \be , g, a_0 , \a )} ( \X _{\cdot } ) \over r + a_0 + X_{\cdot , \nu } + a_{\cdot }} \Big\} \Big] \Big] \non \\
&= E \Big[ \sum_{\nu = 1}^{n} \sum_{i = 1}^{m} \Big[ - \ph _{i, \nu }^{( a_0 , \a )} \phi _{\nu }^{( \al , \be , g, a_0 , \a )} ( \X ) \non \\
&\quad + \ph _{i, \nu }^{\rm{U}} \log \Big\{ 1 + {1 \over r + a_0 + X_{\cdot , \nu } + a_{\cdot } - 1} {K( \al + 1, \be , g, r + a_0 , \X _{\cdot } + a_{\cdot } \j ^{(N)} ) \over K( \al , \be , g, r + a_0 , \X _{\cdot } + a_{\cdot } \j ^{(N)} )} \Big\} \Big] \Big] \non \\
&= E[ - I_{1, n}^{( \al , \be , g, a_0 , \a )} ( \x ) + I_{2, n}^{( \al , \be , g, a_0 , \a )} ( \x ) ] \text{,} \non 
\end{align}
where 
\begin{align}
I_{1, n}^{( \al , \be , g, a_0 , \a )} ( \x ) &= \sum_{\nu = 1}^{n} {x_{\cdot , \nu } + a_{\cdot } \over r + a_0 + x_{\cdot , \nu } + a_{\cdot }} {K( \al + 1, \be , g, r + a_0 , \x _{\cdot } + a_{\cdot } \j ^{(N)} + \e _{\nu }^{(N)} ) \over K( \al , \be , g, r + a_0 , \x _{\cdot } + a_{\cdot } \j ^{(N)} )} \text{,} \non \\
I_{2, n}^{( \al , \be , g, a_0 , \a )} ( \x ) &= \sum_{\nu = 1}^{n} {x_{\cdot , \nu } \over r + x_{\cdot , \nu } - 1} \log \Big\{ 1 + {1 \over r + a_0 + x_{\cdot , \nu } + a_{\cdot } - 1} {K( \al + 1, \be , g, r + a_0 , \x _{\cdot } + a_{\cdot } \j ^{(N)} ) \over K( \al , \be , g, r + a_0 , \x _{\cdot } + a_{\cdot } \j ^{(N)} )} \Big\} \text{,} \non 
\end{align}
and $\x _{\cdot } = ( x_{\cdot , 1} , \dots , x_{\cdot , N} )' = \big( \sum_{i = 1}^{m} x_{i, 1} , \dots , \sum_{i = 1}^{m} x_{i, N} \big) '$ for $\x = ( x_{i, \nu } )_{1 \le i \le m, \, 1 \le \nu \le N} \in {\mathbb{N} _0}^{m \times N}$. 
Since $- I_{1, n}^{( \al , \be , g, a_0 , \a )} ( \bm{0} ^{(m, N)} ) + I_{2, n}^{( \al , \be , g, a_0 , \a )} ( \bm{0} ^{(m, N)} ) < 0$, it is sufficient to show that $- I_{1, n}^{( \al , \be , g, a_0 , \a )} ( \x ) + I_{2, n}^{( \al , \be , g, a_0 , \a )} ( \x ) \le 0$ for all $\x \in {\mathbb{N} _0}^{m \times N} \setminus \{ \bm{0} ^{(m, N)} \} $.

Fix $\x = ( x_{i, \nu } )_{1 \le i \le m, \, 1 \le \nu  \le N} \in {\mathbb{N} _0}^{m \times N} \setminus \{ \bm{0} ^{(m, N)} \} $. 
Let $z _{\nu } = \sum_{i = 1}^{m} x_{i, \nu }$ for $\nu = 1, \dots , N$ and let $\z = ( z_1 , \dots , z_N )'$ and $z = \sum_{\nu = 1}^{N} z_{\nu }$. 
We use the abbreviated notation 
\begin{gather}
\tilde{I} _1 = I_{1, n}^{( \al , \be , g, a_0 , \a )} ( \x ) \text{,} \quad \tilde{I} _2 = I_{2, n}^{( \al , \be , g, a_0 , \a )} ( \x ) \text{,} \quad \tilde{I} = - \tilde{I} _1 + \tilde{I} _2 \text{,} \non \\
\tilde{K} (l) = K( \al + l, \be , g, r + a_0 , \x _{\cdot } + a_{\cdot } \j ^{(N)} ) \text{,} \non \\
\tilde{K} (l, \nu ) = K( \al + l, \be , g, r + a_0 , \x _{\cdot } + a_{\cdot } \j ^{(N)} + \e _{\nu }^{(N)} ) \text{,} \non 
\end{gather}
for $l = 0, 1, 2$ and $\nu = 1, \dots , N$. 
Also, let 
\begin{align}
f_{\al , \be , g, a_0 , \a } (t) &= t^{\al - 1} e^{- \be t} g(t) \prod_{\nu = 1}^{N} {\Ga (t + r + a_0 ) \over \Ga (t + r + a_0 + z_{\nu } + a_{\cdot } )} \non 
\end{align}
for $t \in (0, \infty )$. 

Clearly, 
\begin{align}
\tilde{I} _2 &= \sum_{\nu = 1}^{n} {z_{\nu } \over r + z_{\nu } - 1} \log \Big\{ 1 + {1 \over r + a_0 + z_{\nu } + a_{\cdot } - 1} {\tilde{K} (1) \over \tilde{K} (0)} \Big\} \non \\
&\le \sum_{\nu = 1}^{n} {z_{\nu } \over r + z_{\nu } - 1} {1 \over r + a_0 + z_{\nu } + a_{\cdot } - 1} {\tilde{K} (1) \over \tilde{K} (0)} \non \\
&\le \sum_{\nu = 1}^{n} {z_{\nu } + a_0 + a_{\cdot } + 2 \over (r + a_0 + z_{\nu } + a_{\cdot } )^2} {\tilde{K} (1) \over \tilde{K} (0)} \label{tKLp1} 
\end{align}
by assumption. 
On the other hand, 
\begin{align}
\tilde{I} _1 &= \sum_{\nu = 1}^{n} {z_{\nu } + a_{\cdot } \over r + a_0 + z_{\nu } + a_{\cdot }} {\tilde{K} (1, \nu ) \over \tilde{K} (0)} \non \\
&= \sum_{\nu = 1}^{n} {z_{\nu } + a_{\cdot } \over (r + a_0 + z_{\nu } + a_{\cdot } )^2} {\tilde{K} (1) \over \tilde{K} (0)} - \sum_{\nu = 1}^{n} {z_{\nu } + a_{\cdot } \over (r + a_0 + z_{\nu } + a_{\cdot } )^2} {\tilde{K} (2, \nu ) \over \tilde{K} (0)} \non \\
&\ge \sum_{\nu = 1}^{n} {z_{\nu } + a_{\cdot } \over (r + a_0 + z_{\nu } + a_{\cdot } )^2} {\tilde{K} (1) \over \tilde{K} (0)} - {1 \over n} \sum_{\nu = 1}^{n} {1 \over (r + a_0 + z_{\nu } + a_{\cdot } )^2} \sum_{\nu = 1}^{n} ( z_{\nu } + a_{\cdot } ) {\tilde{K} (2, \nu ) \over \tilde{K} (0)} \label{tKLp2} 
\end{align}
by the covariance inequality. 
Furthermore, by integration by parts, we have 
\begin{align}
( \al + 1) \tilde{K} (1) &= \int_{0}^{\infty } t^2 f_{\al , \be , g, a_0 , \a } (t) \Big\{ \be + {- g' (t) \over g(t)} \Big\} dt \non \\
&\quad + \sum_{\nu = 1}^{N} \int_{0}^{\infty } t^2 f_{\al , \be , g, a_0 , \a } (t) \Big\{ {{\Ga }' (t + r + a_0 + z_{\nu } + a_{\cdot } ) \over \Ga (t + r + a_0 + z_{\nu } + a_{\cdot } )} - {{\Ga }' (t + r + a_0 ) \over \Ga (t + r + a_0 )} \Big\} \non \\
&\ge \sum_{\nu = 1}^{n} ( z_{\nu } + a_{\cdot } ) \tilde{K} (2, \nu ) \text{,} \label{tKLp3} 
\end{align}
where the equality follows since $\Ga (t) \sim t^{- 1}$ as $t \to 0$ while $\prod_{\nu = 1}^{N} \{ \Ga (t + r + a_0 ) / \Ga (t + r + a_0 + z_{\nu } + a_{\cdot } ) \} \sim t^{- z - N a_{\cdot }}$ as $t \to \infty $ and where the inequality follows from Lemma \ref{lem:digamma}. 
Hence, combining (\ref{tKLp1}), (\ref{tKLp2}), and (\ref{tKLp3}), we obtain 
\begin{align}
\tilde{I} &\le - \sum_{\nu = 1}^{n} {z_{\nu } + a_{\cdot } - ( \al + 1) / n \over (r + a_0 + z_{\nu } + a_{\cdot } )^2} {\tilde{K} (1) \over \tilde{K} (0)} + \sum_{\nu = 1}^{n} {z_{\nu } + a_0 + a_{\cdot } + 2 \over (r + a_0 + z_{\nu } + a_{\cdot } )^2} {\tilde{K} (1) \over \tilde{K} (0)} \non \\
&=  - \sum_{\nu = 1}^{n} {- a_0 - 2 - ( \al + 1) / n \over (r + a_0 + z_{\nu } + a_{\cdot } )^2} {\tilde{K} (1) \over \tilde{K} (0)} \text{,} \non 
\end{align}
the right-hand side of which is nonpositive by assumption (\ref{eq:assumption_KL}), and the result follows. 
\hfill$\Box$

\bigskip

\noindent
{\bf Proof of Proposition \ref{prp:MCMC}.} \ \ Properties (ii) and (iv) are trivial. 
Property (iii) follows since 
\begin{align}
\int_{0}^{\infty } \pi ( \p , t | \al , \be , a_0 , \a _1 , \dots , \a _N ) dt = \pi ( \p | \al , \be , a_0 , \a _1 , \dots , \a _N ) \non 
\end{align}
for $\p \in {D_m}^N$. 
For part (i), note that the integrals are finite only if $a_0 \ge 0$ since otherwise 
\begin{align}
\int_{{D_m}^N} \pi ( \p | \al , \be , a_0 , \a _1 , \dots , \a _N ) d\p \ge \int_{{D_m}^N} \pi _{\al , \be , g_1 , a_0 , \overline{a} \j ^{(m)}} ( \p ) d\p = \infty \text{,} \non 
\end{align}
where $\overline{a} = \max \{ \max \{ a_{1, 1} , \dots , a_{m, 1} \} , \dots , \max \{ a_{1, N} , \dots , a_{m, N} \} \} $, by Lemma \ref{lem:propriety}. 
Suppose that $a_0 \ge 0$. 
Then we have 
\begin{align}
&\int_{{D_m}^N \times (0, \infty )} \pi ( \p , t | \al , \be , a_0 , \a _1 , \dots , \a _N ) d( \p , t) / \prod_{\nu = 1}^{N} \prod_{i = 1}^{m} \Ga( a_{i, \nu } ) \non \\
&= \int_{0}^{1} t^{\al - 1} e^{- \be t} \Big\{ \prod_{\nu = 1}^{N} {\Ga (t + a_0 ) \over \Ga (t + a_0 + a_{\cdot , \nu } )} \Big\} dt + \int_{1}^{\infty } t^{\al - 1} e^{- \be t} \Big\{ \prod_{\nu = 1}^{N} {\Ga (t + a_0 ) \over \Ga (t + a_0 + a_{\cdot , \nu } )} \Big\} dt \text{.} \non 
\end{align}
The first term on the right side is finite if and only if $a_0 > 0$ or $\al > N$ since $\Ga (t) \sim t^{- 1}$ as $t \to 0$. 
The second term on the right side is finite if and only if $a_{\cdot , \cdot } > \al $ or $\be > 0$ since $\prod_{\nu = 1}^{N} \Ga (t + a_0 ) / \Ga (t + a_0 + a_{\cdot , \nu } ) \sim t^{- a_{\cdot , \cdot }}$ as $t \to \infty $. 
This completes the proof. 
\hfill$\Box$

\section*{Acknowledgments}
Research of the second author was supported in part by Grant-in-Aid for Scientific Research (18K11188, 15H01943 and 26330036) from Japan Society for the Promotion of Science.

\end{document}